\RequirePackage{fix-cm}
\documentclass[review, twocolumn]{svjour3}          
\smartqed  
\usepackage{graphicx}

\usepackage{lineno,hyperref}
\modulolinenumbers[5]

\usepackage{algorithm}
\usepackage{algorithmic}
\usepackage{cite}
\usepackage{listings}
\usepackage{subfigure}
\usepackage{graphicx}
\usepackage{array}
\usepackage{multirow}
\newcommand{\specialcell}[2][l]{\begin{tabular}[#1]{@{}#1@{}}#2\end{tabular}}
\usepackage{amsmath}

\newcommand{\of}[1]{\left( #1 \right)}

\newcommand{\eps}{\mbox{$\epsilon$}}

\let\vec=\mathbi%
\let\mat=\mathbf%
\let\set= \mathcal%

\usepackage{mleftright}
\usepackage{amssymb}
\usepackage{color}
\usepackage{tikz}
\usetikzlibrary{matrix}

\definecolor{RED}{rgb}{1,0,0}
\definecolor{GREEN}{rgb}{0.8, 0, 0.5}
\definecolor{BLUE}{rgb}{0,0.4,1}

\newcommand\given[1][]{\:#1\vert\:}

\DeclareMathOperator*{\argmax}{arg\,max}

\journalname{Computing and Visualization in Science}

\begin{document}

\title{Sobol Tensor Trains for Global Sensitivity Analysis
\thanks{This work was partially supported by the UZH Forschungskredit ``Candoc'', grant number FK-16-012.}
}

\author{Rafael Ballester-Ripoll \and
Enrique G. Paredes \and
Renato Pajarola}


\institute{R. Ballester-Ripoll, E.G. Paredes, and R. Pajarola \at
              Visualization and MultiMedia Lab, Dept. of Informatics,\\ University of Zurich. Switzerland\\
              \email{rballester@ifi.uzh.ch, egparedes@ifi.uzh.ch,\\ pajarola@ifi.uzh.ch}         
}

\date{Received: date / Accepted: date}

\maketitle

\begin{abstract}
Sobol indices are a widespread quantitative measure for variance-based global sensitivity analysis, but computing and utilizing them remains challenging for high-dimensional systems. 
We propose the tensor train decomposition (TT) as a unified framework for surrogate modeling and global sensitivity analysis via Sobol indices. We first overview several strategies to build a TT surrogate of the unknown true model using either an adaptive sampling strategy or a predefined set of samples.
We then introduce and derive the Sobol tensor train, which compactly represents the Sobol indices for all possible joint variable interactions which are infeasible to compute and store explicitly.
Our formulation allows efficient aggregation and subselection operations: we are able to obtain related indices (closed, total, and superset indices) at negligible cost. 
Furthermore, we exploit an existing global optimization procedure within the TT framework for variable selection and model analysis tasks. We demonstrate our algorithms with two analytical engineering models and a parallel computing simulation data set.

\end{abstract}

\section{Introduction}

A crucial task when analyzing computational models and physical simulations is assessing the influence of each input variable (and all combinations thereof) on the model's output. The quantitative study of such influences is known as \emph{sensitivity analysis} (SA). When the variables themselves are stochastic, the propagation of their uncertainty towards the model output must also be taken into account. We focus on variance-based SA, often referred to as \emph{analysis of variances} (ANOVA), and in particular the so-called \emph{Sobol decomposition}. It approximates the parametrized model as a sum of simpler functions, each depending on only a subset of the original set of variables. The sensitivity to each variable is then reflected by the functions that depend on it, and can therefore be estimated as their relative contribution to the output's overall statistical variance. These relative variances are known as \emph{Sobol indices} and have become a standard tool for global SA in the last few decades~\cite{SRACCGST:08, Sudret:08, MILR:09, ODC:14, IL:15}.

A popular method to compute such indices is via \emph{Monte Carlo} (MC) integration estimators on a suitable set of samples within the variable space (the \emph{sampling plan}). This was already outlined in the original paper by Sobol~\cite{Sobol:90} and gained momentum thereafter. However, MC convergence is slow w.r.t the number of samples available~\cite{IL:15}. Structured sampling plans exist that improve convergence, e.g. \emph{Latin hypercube sampling} or \emph{quasi-random sequences} (quasi-MC). If needed one may favor estimators for \emph{total effect} indices, i.e. \emph{quantities of interest} (QoI) that aggregate Sobol indices of diverse orders together. Unfortunately, a plan tailored to estimate a particular index, or set thereof, may be suboptimal for other indices. In practice, analysts often choose to restrict the Sobol decomposition to interactions of low-order (e.g. up to 2), and/or perform a prior dimensionality reduction in what is known as \emph{screening} (e.g. freezing seemingly unimportant variables). Such simplifications greatly reduce the computational complexity, but pose a risk: they might fail to detect significant complex interactions between variables, and over-zealous reduction can harm subsequent processing steps in the analysis pipeline.

A complementary approach to direct MC estimation is building a \emph{surrogate model}, also known as \emph{response surface model} or \emph{metamodel}, in an offline preliminary step. The surrogate acts as an interpolator that is fast to evaluate and can approximate the true unknown model at arbitrary sampling points~\cite{QHSGVT:05}. This strategy is attractive when sample acquisition is expensive or highly dynamic, especially if the analyst would like to estimate new indices on demand. Furthermore, several surrogates can produce Sobol indices in a more direct manner~\cite{IL:15}, thus avoiding MC integration altogether. However, dealing with high-dimensional parametric systems, i.e. with a significant number $N$ of input variables, remains a major challenge. Even if the chosen surrogate scales well with the dimensionality~\cite{KS:16}, the sheer number of sensitivity indices is by definition exponential, as there are $2^N-1$ Sobol indices, out of which $\binom{N}{M}$ for any fixed order $1 \le M \le N$ may be chosen. For moderate or large values of $N$, general queries of the form \emph{``retrieve the largest indices of any order''} or \emph{``compute the total variance for interactions of order up to $k$''}
quickly become computationally intractable.

To address these problems we propose a data structure that compactly stores \emph{all} Sobol sensitivity indices in a compressed form. It is based on tensor decompositions, in particular the tensor train (TT)~\cite{Oseledets:11} method. The TT format is designed to avoid the curse of dimensionality and excels in high-dimensional approximation and compression in several fields ranging from physics and quantum chemistry to engineering and data mining. It can be built using various sampling/interpolation techniques, both when a sample acquisition plan is required and when the set of known samples is fixed beforehand.

\begin{figure*}[t]\centering
  \includegraphics[width=2\columnwidth]{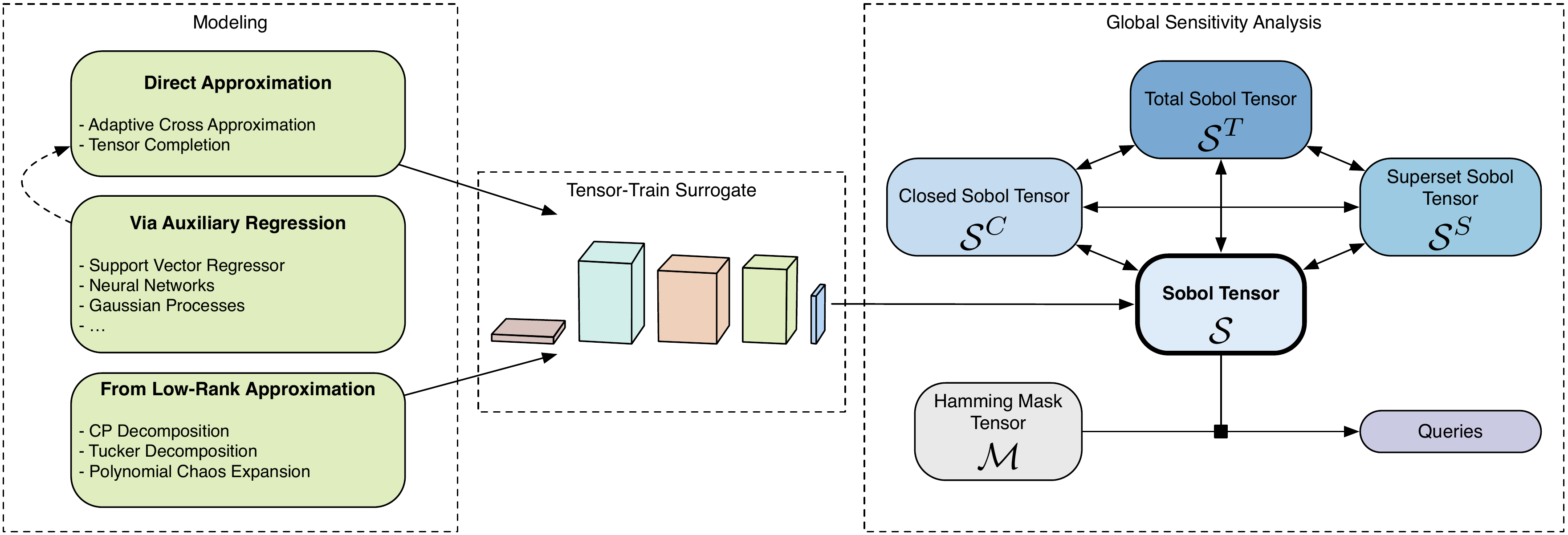}
  \caption{Pipeline for TT-based global sensitivity analysis: a model with $N$ input variables is approximated as an $N$-dimensional tensor, from which we extract a compact $2^N$ tensor $\set{S}$ approximating all $2^N-1$ Sobol indices. This tensor can be then used for various aggregation, analysis and query/optimization tasks.}
  \label{fig:sensitivity_analysis}
\end{figure*}

While formulas to compute individual or aggregated Sobol indices from various low-rank surrogates have been already derived in the recent literature~\cite{Rai:14, DKLM:14, KS:16}, our approach is the first to assemble the complete set of all indices in a unified and compact tensor format that can be manipulated and queried for statistics, model reduction, visual analytics, and more. See Fig.~\ref{fig:sensitivity_analysis} for a summary diagram of our pipeline.


\subsection{Contributions}

We introduce the Sobol tensor $\set{S}$, an $N$-dimensional TT-compressed multiarray encoding all possible $2^N-1$ Sobol indices for global SA, and show its derivation from an arbitrary TT surrogate model. We further extract the related aggregated tensors $\set{S}^S$, $\set{S}^C$ and $\set{S}^T$ containing the \emph{superset}, \emph{closed} and \emph{total} indices, respectively. All these indices can be derived from each other via union/intersection operations that are translated to the tensor-compressed domain as simple matrix additions and subtractions.
By combining these ideas with existing optimization algorithms for the TT format we are able to answer several computationally challenging types of global SA queries that often arise during variable selection and model interpretation.



\subsection{Notation}

Multidimensional arrays, herein called \emph{tensors}, have sizes denoted by $I_1, ..., I_N$ where $N$ is the number of dimensions. \emph{Tensor ranks} generalize the matrix rank for $N > 2$ and reflect, in a sense, the complexity of a tensor; they use the symbols $R_n$. For simplicity we sometimes use $I := \max \{I_n\}_n$ and $R := \max \{R_n\}_n$. Vectors, matrices and tensors use bold lowercase, bold uppercase and calligraphic letters as in $\vec{x}$, $\mat{U}$ and $\set{T}$ respectively. Their elements are denoted by square brackets with indices starting from 0, following NumPy convention (e.g. $\mat{U}[1, 0]$). Furthermore, we refer to the element-wise (i.e., Hadamard) product of two tensors as $\set{A} \circ \set{B}$, whereas the Kronecker product of matrices is written as $\mat{A} \otimes \mat{B}$. We use the Frobenius norm $\| \cdot \| := \| \cdot \|_2$, and accuracy between groundtruth data $\set{A}$ and an approximation $\set{B}$ is measured with the \emph{relative error}: $\eps_{\set{A}}(\set{B}) := \|\set{A} - \set{B}\| / \|\set{A}\| \ge 0$.

We denote tuples of indices as $\pmb{\alpha} \in \{0, 1\}^N$ and $\pmb{\alpha} \subseteq \{1, ..., N\}$ interchangeably: a 0 (resp. 1) in the former notation means an index is absent (resp. present) in the latter. If a function $f: \mathbb{R}^4 \to \mathbb{R}$ only depends \emph{effectively} on the two last variables, we may alternatively write $\pmb{\alpha} = [0, 0, 1, 1]$ or $\pmb{\alpha} = \{3, 4\}$, and $f_{\pmb{\alpha}}(\vec{x}) \equiv f_{\pmb{\alpha}}(\vec{x}_{\pmb{\alpha}}) \equiv f_{3,4}(x_3, x_4)$ similarly to~\cite{Sudret:08, Owen:14}. Cardinality of a set of variables is denoted as $|\pmb{\alpha}|$, and coincides with the \emph{Hamming weight} (bit sum) of its binary representation. Last, we write tuple complements as $-\pmb{\alpha} \equiv \{1, ..., N\} \setminus \pmb{\alpha}$.

\section{Variance-Based Sensitivity Analysis} \label{sec:related_work}

\subsection{Sobol Decomposition} \label{sec:sobol_decomposition}

Variance-based SA dates back to the early 20th century and comprises a set of related techniques for statistical analysis of multidimensional data, out of which the ANOVA is arguably the most widely known. The functional ANOVA decomposition~\cite{ES:81}, also known as high-dimensional model representation~\cite{BK:15} or Sobol decomposition~\cite{Sobol:90}, writes any integrable multidimensional function over a rectangle $f: \Omega = \Omega_1 \times ... \times \Omega_N \subset \mathbb{R}^{N} \rightarrow \mathbb{R}$ as the following sum of subfunctions:

\begin{equation}
\label{eq:anova1}
\begin{aligned}
	f(\vec{x}) = f_{\emptyset} + \sum_{i=1}^N f_i(x_i) + \sum_{1 \le i < j \le N} f_{ij}(x_i, x_j) + \dots \\ + \sum_{\substack{\pmb{\alpha} \subset \{1,...,N\} \\ |\pmb{\alpha}| = n}} f_{\pmb{\alpha}}(x_{\pmb{\alpha})} + \dots + f_{1,...,N}(x_1,...,X_N)
\end{aligned}
\end{equation}
where each $f_{\pmb{\alpha}}$ only depends effectively on the indices contained in $\pmb{\alpha}$. The $f_{\pmb{\alpha}}$ are uniquely determined if orthogonality w.r.t. a separable measure $dF(\vec{x}) = dF_1(x_1) \\ \cdots dF_N(x_N)$ is imposed:

\begin{equation}
\int_{\Omega} f_{\pmb{\alpha}}(\vec{x}_{\pmb{\alpha}}) f_{\pmb{\beta}}(\vec{x}_{\pmb{\beta}}) \, dF(\vec{x}) = 0 \mbox{ for any } \pmb{\alpha} \ne \pmb{\beta}
\end{equation}
Then, the explicit decomposition stems from iterative integrations and subtractions~\cite{Sobol:90}:

\begin{equation}
\label{eq:anova2}
\begin{split}
f_{\emptyset} = \int_{\Omega} f(\vec{x}) \, dF(\vec{x}) \\
f_n(x_n) = \int_{\Omega_{-n}} f(\vec{x}) \, dF_{-n}(\vec{x}_{-n}) - f_{\emptyset} \\
f_{nm}(x_n, x_m) = \int_{\Omega_{-nm}} f(\vec{x}) \, dF_{-nm}(\vec{x}_{-nm}) \\ - f_n(x_n) - f_m(x_m) - f_{\emptyset} \\
\cdots \\
f_{\pmb{\alpha}}(\vec{x}_{\pmb{\alpha}}) = \int_{\Omega-{\pmb{\alpha}}} f(\vec{x}) \, dF_{-\pmb{\alpha}}(\vec{x}_{-\pmb{\alpha}}) - \sum_{\pmb{\beta} | \pmb{\beta} \subset \pmb{\alpha}} f_{\pmb{\beta}}(\vec{x}_{\pmb{\beta}})
\end{split}
\end{equation}

Eqs.~\ref{eq:anova1} to ~\ref{eq:anova2} are useful in the context of uncertainty quantification, namely when one has a model depending on $N$ independent random variables $x_1, ..., x_N$. Under this assumption, their joint \emph{probability density function} (PDF) plays the role of our separable measure in $\Omega$ and the integrals are expectations of each subfunction, starting with $f_{\emptyset} = \operatorname{E}[f]$. Eq.~\ref{eq:anova2} is then

\begin{equation} \label{eq:expected_values}
f_{\pmb{\alpha}}(\vec{x}_{\pmb{\alpha}}) = (\operatorname{E}_{-\pmb{\alpha}} [f])(\vec{x}_{\pmb{\alpha}}) - \sum_{\pmb{\beta} | \pmb{\beta} \subset \pmb{\alpha}} f_{\pmb{\beta}}(\vec{x}_{\pmb{\beta}})
\end{equation}

\subsection{Variance and Sobol Indices} \label{sec:variance_and_sobol}

The \emph{variance indices} $D_{\pmb{\alpha}}$ are defined as the variance contributed by each of the $f_{\pmb{\alpha}}$, w.r.t. the measure $F$. Thus, the Sobol decomposition builds up a partition of the overall variance $D$:

\begin{equation}
D = \sum_{\pmb{\alpha}} D_{\pmb{\alpha}} = \operatorname{V}[f] = \operatorname{E}[f^2] - \operatorname{E}^2[f] = \int_{\Omega} f(\vec{x})^2 dF(\vec{x}) - f_{\emptyset}^2
\end{equation}

The \emph{Sobol indices}~\cite{Sobol:90} in turn map the relative variances w.r.t. the total model variance:

\begin{equation}
\begin{array}{c}
S: \mathcal{P}(\{1, ..., N\}) \setminus \emptyset \to [0, 1] \\
S_{\pmb{\alpha}} := D_{\pmb{\alpha}}/D \\
\sum_{\pmb{\alpha}} S_{\pmb{\alpha}} = 1
\end{array}
\end{equation}

These indices are an invaluable tool in many SA settings~\cite{STCR:04}, for example in factor prioritization (reducing uncertainty), factor fixing (identifying non-influential variables), risk minimization, reliability engineering, etc. They are also helpful to select good dimension orderings that lead to more compact surrogate models (example 5.8 by Bigoni~\cite{Bigoni:14}; also considered in~\cite{DKLM:14}). They are hyperedges of a hypergraph, since they encode $n$-ary relations within subsets of $\{1, ..., N\}$. Alternatively they can be thought of in terms of set cardinalities, as the sum of all $S_{\pmb{\alpha}}$ equals $1$ (see e.g.~\cite{Owen:13} and \cite{SRACCGST:08}, Sec.~1.2.15).

Several surrogate models lend themselves well to direct estimation of Sobol indices. Examples in the literature include PCE of bounded degree~\cite{Sudret:08}, low-rank sums of separable PCE terms~\cite{KS:16}, Gaussian processes~\cite{MILR:09}, TT~\cite{DKLM:14}, spectral TT~\cite{BEM:16}, etc. However, there are $2^N - 1$ possible QoI after excluding the trivial tuple $\pmb{\alpha} = [0, ..., 0] \equiv \emptyset$. As $N$ grows, this magnitude poses challenges in both the computational and the model interpretation aspects.


\subsection{Related Indices} \label{sec:related_indices}

One may derive alternative indices by adding and/or subtracting together the standard $S_{\pmb{\alpha}}$, effectively configuring a set algebra.

\subsubsection{Total Indices}

Denoted as $S^T_{\pmb{\alpha}}$, they are also called \emph{upper indices}~\cite{Owen:13}. They represent all joint indices that include any variable from $\pmb{\alpha}$:

\begin{equation}
S^T_{\pmb{\alpha}} := \sum_{\pmb{\beta} | \pmb{\alpha} \cap \pmb{\beta} \ne \emptyset} S_{\pmb{\beta}}
\end{equation}

For example, in a 3-variable model we have $S^T_{1,2} = S_1 + S_2 + S_{1,2} + S_{1,3} + S_{2,3} + S_{1,2,3}$.
If $|\pmb{\alpha}| = 1$ we are encoding the total influence of a single variable also accounting for its higher-order interactions with all other variables. In this case the indices are called sometimes \emph{total effects}~\cite{HS:96}, and have been used to identify and select the most relevant variables, for example by sorting $S_n^T$ and choosing the $k$ largest~\cite{Fock:14, AJKS:17}. However, this criterion may lead to overestimating variables that exhibit large overlapping variance contributions. 

\subsubsection{Closed Indices}

Denoted as $S^C_{\pmb{\alpha}}$, they are also called \emph{first-order indices}~\cite{Sudret:08} or \emph{lower indices}~\cite{Owen:13}. They sum the variance contributions of all non-empty tuples contained in $\pmb{\alpha}$:

\begin{equation}
S^C_{\pmb{\alpha}} := \sum_{\pmb{\beta} | \pmb{\alpha} \supseteq \pmb{\beta}} S_{\pmb{\beta}}
\end{equation}

For example, for 3 variables we have $S^C_{1,2} = S_1 + S_2 + S_{1,2}$. Also, for any single variable $n$ we have $S^C_n = S_n$. The closed indices can be written in terms of the total indices as $S^C_{\pmb{\alpha}} = 1 - S^T_{-\pmb{\alpha}}$.

\subsubsection{Superset Indices}

The $\set{S}^S$ aggregate all indices that contain a tuple~\cite{Hooker:04}:

\begin{equation}
S^S_{\pmb{\alpha}} := \sum_{\pmb{\beta} | \pmb{\alpha} \subseteq \pmb{\beta}} S_{\pmb{\beta}}
\end{equation}

For example, $S^S_{1,2} = S_{1,2} + S_{1,2,3}$.





\section{Tensor Approximation}
\label{sec:tensor_decompositions}


Decomposing multidimensional arrays (tensors) in terms of simpler, separable terms is a fruitful approach in compression, interpolation and metamodeling applications, and their fundamentals reach out to other important mathematical frameworks including principal component analysis, wavelet transforms, polynomial chaos, etc. We briefly introduce first CP and Tucker since they are arguably the two most popular tensor models, and we support them in Sec.~\ref{sec:tt_construction} as optional intermediate surrogate models. We cover then the more recent TT decomposition, which is the keystone of all algorithms presented in this paper. The section concludes with related work on tensor-based surrogate modeling and SA.

\subsection{CANDECOMP/PARAFAC}

The CP, also known as \emph{canonical} or \emph{Kruskal decomposition}, is the earliest and most straightforward extension of the singular value decomposition (SVD) for more than 2 dimensions~\cite{Harshman:70}. It approximates a tensor $\set{T}$ element-wise as follows:

\begin{equation}
\label{eq:cp}
\set{T}[x_1, ..., x_N] \approx \sum_{r=1}^R \lambda_r \cdot \mat{U}^{(1)}[x_1, r] \cdot ... \cdot \mat{U}^{(N)}[x_N, r]
\end{equation}
where $R$ is the \emph{CP rank}. The $\mat{U}^{(n)}$ are known as \emph{factor matrices} or simply \emph{factors}. We can write Eq.~\ref{eq:cp} more compactly using double bracket notation as $\set{T} \approx [[\pmb{\lambda}; \mat{U}^{(1)}, ..., \mat{U}^{(N)}]]$. The $\pmb{\lambda}$ can be optionally absorbed column-wise by the factors $\mat{U}^{(n)}$ and omitted in the notation. Unfortunately, the set of $N$-dimensional tensors of fixed CP-rank $R$ is not closed in $\mathbb{R}^N$, and finding the best rank-$R$ approximation of a given tensor is an ill-posed problem~\cite{SL:08}. On the positive side, the CP format needs $O(I N R)$ elements for storage, i.e. it is linear w.r.t. the number of dimensions.

\subsection{Tucker}

The Tucker decomposition~\cite{Tucker:66, LMV:00b} is also known as \emph{higher-order SVD} (HOSVD), \emph{N-mode PCA}, and \emph{low multilinear rank approximation} (LMLRA). It extends CP by considering all interactions between its factor columns, weighted by an $N$-dimensional core $\set{B}$ of size $R_1 \times ... \times R_N$:

\begin{equation}
\set{T}[x_1, ..., x_N] \approx \sum_{\vec{r}=\vec{1}}^{\vec{R}} \set{B}[\vec{r}] \cdot \mat{U}^{(1)}[x_1, r_1] \cdot ... \cdot \mat{U}^{(N)}[x_N, r_N]
\end{equation}

Approximating a tensor with the Tucker format is a stable procedure~\cite{LMV:00b}. However, $O(R^N + I N R)$ elements need to be stored, i.e. it still suffers from the curse of dimensionality. It is therefore mostly used for up to a handful of dimensions only. The Tucker model is related to \emph{polynomial chaos expansions} (PCE), as we detail later in Sec.~\ref{sec:from_pce}.

\subsection{Tensor Train}

\begin{figure}[ht]\centering 
  \includegraphics[width=0.8\columnwidth]{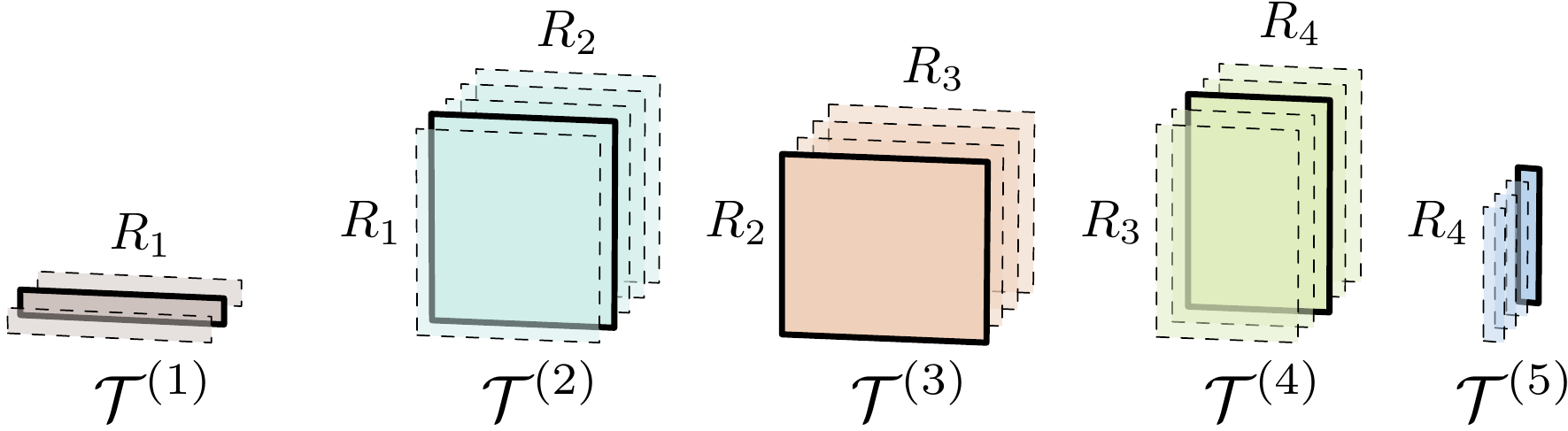}
  \caption{A 5D tensor train of size $3 \times 5 \times 4 \times 5 \times 4$. By multiplying the highlighted matrices together we obtain the element $\set{T}^{(1)}[1] \cdot \set{T}^{(2)}[1] \cdot \set{T}^{(3)}[0] \cdot \set{T}^{(4)}[2] \cdot \set{T}^{(5)}[4] = \set{T}[1, 1, 0, 2, 4]$.}
  \label{fig:tt}
\end{figure}

The TT model (also known since the 1990s as \emph{matrix product states} or \emph{linear tensor network}) was recently rediscovered by Oseledets~\cite{Oseledets:11}. It aims to unite the advantages of both CP and Tucker, especially for high $N$. It uses a sequence of 3D cores, compactly written as $[[\set{T}^{(1)},\dots,\set{T}^{(n)}]]$. The reconstruction is a sequence of matrix products:

\begin{equation}
\set{T}[x_1,\dots,x_N] = \set{T}^{(1)}[x_1] \cdot ... \cdot \set{T}^{(N)}[x_N]
\end{equation}
where $\set{T}^{(n)}[x_n]$ is a shorthand for the $x_n$-th slice along mode 2, i.e. $\set{T}^{(n)}[:, x_n, :]$. The core dimensions are $R_{n-1} \times I_n \times R_n$ for $n = 1, ..., N$, with $R_0 := R_N := 1$ by convention. The $R_n$ are called \emph{TT ranks} and are bounded from above by CP's $R$ rank. Fig.~\ref{fig:tt} illustrates the structure and indexing of the core sequence.

In contrast to the Tucker model, the TT format needs $O(I N R^2)$ elements and thus grows linearly w.r.t the number of dimensions $N$.




\subsubsection{Operations in the TT Format}

Multiplication by a scalar and tensor-wise addition/product may be achieved by simple manipulations of the TT cores as shown in~\cite{Oseledets:11}; see also~\ref{app:operations} for more in-depth details. Furthermore, thanks to the so-called \emph{adaptive cross-approximation} (ACA) framework in the TT format~\cite{OT:10, SO:11}, these and many other operations can be accomplished in $O(I N R^3)$ operations at most, i.e. devoid of the curse of dimensionality. These include arbitrary element-wise functions, differentiation, integration, convolution, and more~\cite{CLOPZM:16, LC:17}. The ranks may grow as a result of such operations. It is crucial to keep them reasonably low at all stages of any computational pipeline, otherwise the benefits of tensor compression vanish. An error-bounded rounding algorithm called TT-round~\cite{Oseledets:11} exists to re-compress down any tensor when needed.

\subsubsection{Global Optimization} \label{sec:global_optimization}

ACA has been successfully used to find the (approximately) maximal element in modulus of a TT tensor~\cite{DKLM:14, OZTSKS:15}, as it was empirically found that the subtensors accessed during cross-approximation very often contain such maximal elements. The algorithmic variant known as \emph{rectangular maxvol} is a tool even more effective for this task~\cite{MO:15} and is the one we use (as released in~\cite{ttpy}).



\subsection{Tensor Surrogates and Sensitivity Analysis}



Tensor decompositions make for attractive surrogates owing to their natural multidimensionality and fast decompression. Several examples~\cite{EHLMZ:13, LME:13, DKLM:14} target solutions of multiparametric partial differential equations (PDEs). Konakli and Sudret~\cite{KS:16} propose an interpolator via sums of separable PCE-based functions (low-rank approximations, LRA) and showed how to extract Sobol indices out of them. This is related to the CP decomposition, with the main difference that their factors are continuous and are sought within the subspace spanned by a few leading orthogonal polynomials.
Vervliet et al.~\cite{VDSL:14} demonstrate CP-based tensor completion and visualization for the melting point of an alloy, depending on the concentration of its 10 different constituent materials. Ballester-Ripoll et al.~\cite{BPP:16} propose visualization diagrams for TT-format surrogates of several mechanical simulations, emphasizing multidimensionality and real-time reconstruction. 

A few papers have extracted Sobol indices from TT surrogates. Dolgov et al.~\cite{DKLM:14} build their decomposition using ACA and derive properties and statistics including means, covariances, level sets, and individual Sobol indices. Zhang et al.~\cite{ZYOKD:15} developed a hierarchical uncertainty quantification algorithm using TT and PCE to estimate a circuit's response depending on its subcomponents' behavior. Rai~\cite{Rai:14} gives formulas to compute Sobol indices from a range of low-rank approximation surrogates, including TT-based.




\section{Construction of TT Surrogates} \label{sec:tt_construction}

Surrogate-based sensitivity analysis methods are only as good as the approximant's accuracy w.r.t. to the true unknown model. A key part of our pipeline is thus obtaining a high-quality TT interpolant. Fortunately, many models can be accurately represented by a low-rank TT model. For example, multiplicative functions (i.e. with the form $f_1(x_1) \cdots f_N(x_N)$) have exactly rank 1, while additive terms (i.e. $f_1(x_1) + ... + f_N(x_N)$) have exactly rank 2. More generally, we can build TT surrogates in a wide range of settings.


\subsection{Preliminaries: Variable Range Discretization} \label{sec:discretization}

The methods we present are applicable to both continuous and categorical variables, and these two kinds may coexist within one model. However, in order to build the tensor product grid $I_1 \times ... \times I_N$ for our variable space, our TT surrogate $\tilde{f}(\vec{x}) \approx f(\vec{x})$ needs to discretize each continuous variable's domain as a finite set $x_n(1) < ... < x_n(I_n)$. To record or evaluate an entry $\vec{x}$, each coordinate $x_n$ must be first quantized to match the corresponding axis discretization. This is not a problem in practice, and discretizing the variable space is indeed a usual feature of several sampling strategies such as factorial design, Morris' method, one-at-a-time design, etc.~\cite{IL:15}. If needed, the grid can be refined by simply increasing the sampling resolution before building the surrogate, and all important TT operations are linear w.r.t. the spatial dimensions $I_n$. 
For simplicity we use nearest-neighbor interpolation to convert an arbitrary $\vec{x}$ to integer tensor indices $1, ..., I_n$.

\subsection{Construction From a Black-Box System} \label{sec:adaptive_sampling}


\emph{Black-box adaptive sampling} is the scenario in which new samples $\mat{X} = \{\vec{x}_1, ..., \vec{x}_P\}$ are to be chosen and evaluated from scratch with no prior information on the inner workings of the true model. One has freedom to select the set of samples $\mat{X}$, and can do so adaptively in order to minimize the model's generalization error. 
Adaptive cross approximation (ACA) builds a progressive sampling plan on the low-rank assumption; it is an example of design of experiments (DoE). ACA has been an active recent research topic~\cite{OST:08, CC:10, OT:10, SO:11} and has become a key tool to create and manipulate tensors. Recent techniques generalize the maximum-volume (\emph{maxvol} for short) algorithm, which approximates a matrix in terms of a carefully selected subset of its rows and columns. In higher dimensions, ACA constructs the plan by progressively sampling certain \emph{tensor fibers}: sets of samples obtained by fixing all parameters but one. This is a case of one-at-a-time sampling, which can improve the DoE's overall efficiency (see also~\cite{SRACCGST:08}, 2.4.2). Like Latin hypercube sampling, this guarantees that all possible discretized values for every variable are used at least once. In this paper we use an \emph{alternating minimal energy} method to select the fibers~\cite{DS:14}, an algorithm whose implementation has publicly been released as part of the Python ttpy toolbox~\cite{ttpy}.

\subsection{From Categorical Data} \label{sec:categorical_data}

Tensors are discrete data structures indexed by discrete axes and thus support categorical variables in a natural way (consider for example the 2D case: the rank of a matrix is not affected if we permute its columns or rows). Populating the missing entries of a tensor without any prior assumption about smoothness is known as \emph{tensor completion} and is a very convenient tool for regression on categorical variables. It is similar to the better known problem of low-rank matrix completion for $N = 2$, but specific algorithms for $N \ge 3$ of course depend heavily on the particular decomposition format chosen (CP, Tucker, TT, etc.). We have implemented an \emph{alternating least squares} (ALS) completion algorithm in the TT format~\cite{GKT:13} and used it to learn a 14D categorical data set (see Sec.~\ref{sec:gemm}). The algorithm is a form of block coordinate descent, whereby one TT core is optimized at a time and the relative error is provably non-increasing.

\subsection{From an Auxiliary Regressor} \label{sec:auxiliary_regression}

More generally, one may want to interpolate the given training set first with a preferred regression method: support vector machines, radial basis functions, Gaussian processes, etc. One can then approximately transform this auxiliary surrogate into the TT representation by an ACA algorithm. This is a very general approach and is feasible as long as the intermediate regression's output can be approximated well by a surrogate of low TT ranks. Under this assumption, ACA works as a universal tool to reduce any model into the TT format, usable whenever an ad-hoc conversion in the compressed domain (such as the ones discussed next) is not available.


\subsection{From Another Low-Rank Decomposition} \label{sec:from_low_rank}

Several well-known surrogate models are actually based on a low-rank expression, or can easily be cast as a low-rank format. We can convert from these cases more directly instead of relying on the general ACA as just discussed in Sec.~\ref{sec:auxiliary_regression}.


\subsubsection{From CP} \label{sec:from_cp}


TT ranks are bounded from above by CP ranks~\cite{Oseledets:11}, and the proof is constructive: Given a rank-$R$ CP decomposition $[[\mat{U}^{(1)}, ..., \mat{U}^{(N)}]]$, an equivalent TT expression $[[\set{T}^{(1)}, ..., \set{T}^{(N)}]]$ can be straightforwardly built as:

\begin{equation}
\set{T}^{(n)}[x_k] := \mbox{diag}(\mat{U}^{(n)}[x_k, :]), k = 1, ..., I_n
\end{equation}
where the $n$-th core has size $R \times I_n \times R$. Thanks to this we can convert an arbitrary low-rank CP surrogate into our preferred standard TT representation.

\subsubsection{From Tucker} \label{sec:from_tucker}

A TT can be also obtained from a Tucker decomposition, although TT-ranks are not bounded by Tucker ranks (and vice versa). To do this conversion we start with the Tucker approximation formula

\begin{equation}
\tilde{f} \approx \set{T} = [[\set{B}; \mat{U}^{(1)}, ..., \mat{U}^{(N)}]]
\end{equation}
and then compress its core in the TT format:

\begin{equation}
\label{eq:tt_tucker}
\set{T} \approx [[\,[[\set{B}^{(1)}, ..., \set{B}^{(N)}]]; \mat{U}^{(1)}, ..., \mat{U}^{(N)}]]
\end{equation}

By multilinearity, the right-hand side of Eq.~\ref{eq:tt_tucker} equals the following expression

\begin{equation}
\label{eq:ett}
[[\set{B}^{(1)} \times_2 \mat{U}^{(1)}, ..., \set{B}^{(N)} \times_2 \mat{U}^{(N)}]]
\end{equation}

\begin{figure}[ht]\centering
  \includegraphics[width=0.75\columnwidth]{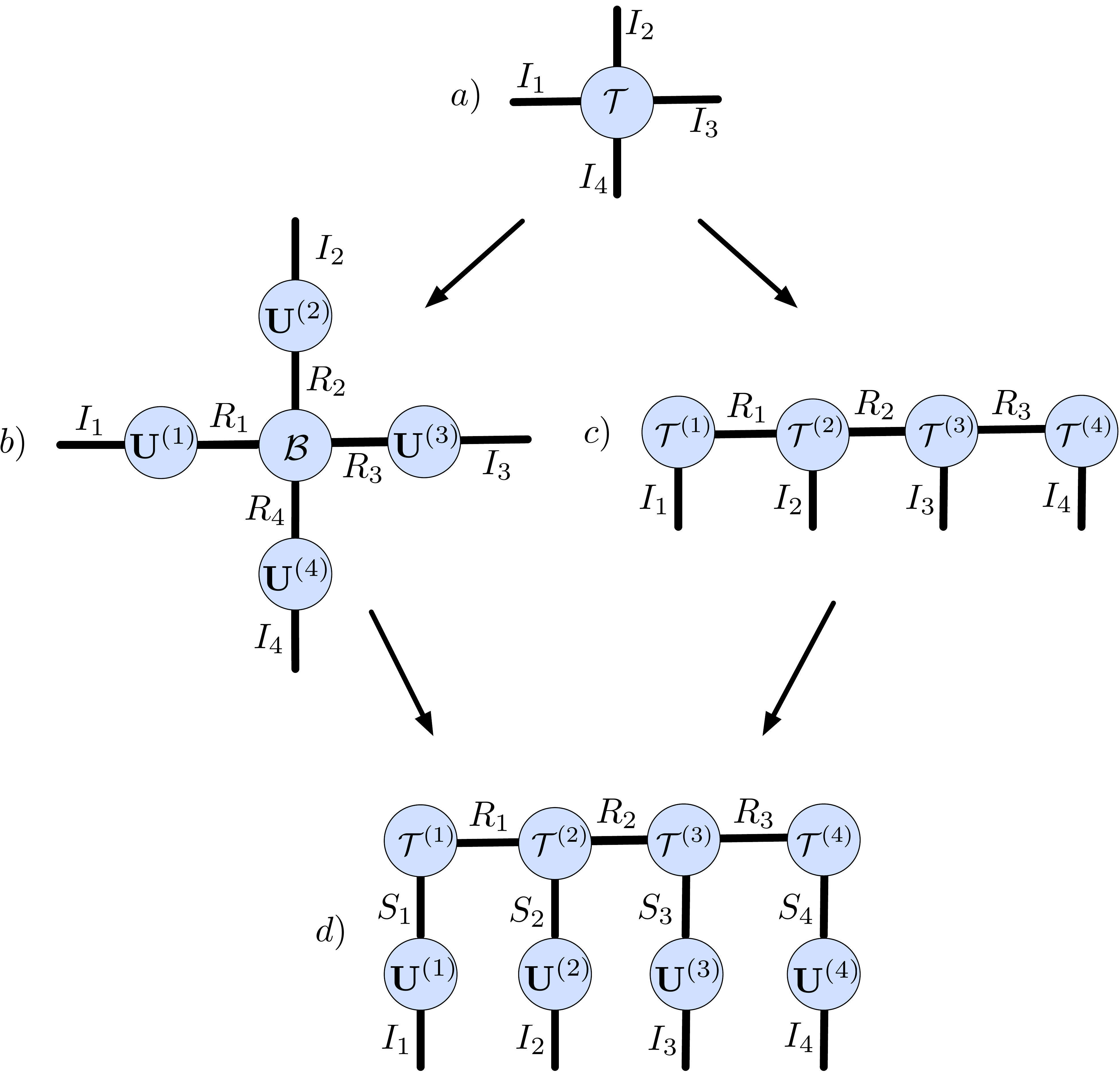}
  \caption{a) A full 4D tensor $\set{T}$; b) in Tucker format; c) in TT format; d) in extended TT format (ETT). The ETT representation can be computed in two alternative but equivalent ways: either via TT compression of the Tucker core $\set{B}$ (left), or via Tucker compression (along the 2nd mode) of each individual TT core $\mathcal{T}^{(n)}$ (right). Similarly, the Tucker format may be cast to TT by following either the path b)-a)-c) (full decompression and compression) or the much less expensive b)-d)-c).}
  \label{fig:commute}
\end{figure}

Eq.~\ref{eq:ett} is a so-called \emph{extended tensor-train} decomposition (ETT), originally defined in~\cite{OT:10b}. See also Fig.~\ref{fig:commute} for a graphical representation in terms of tensor networks, similarly to~\cite{CLOPZM:16}.

The final TT cores are retrieved by explicitly performing the tensor-times-matrix operations:

\begin{equation}
\set{T}^{(n)} = \set{B}^{(n)} \times_2 \mat{U}^{(n)}
\end{equation}
which increases the overall size, but is still linear w.r.t. $N$.

\subsubsection{From Polynomial Chaos Expansions} \label{sec:from_pce}


PCE surrogate models have been used in stochastic modeling and uncertainty quantification for decades. A PCE is based on a set of $N$ polynomial bases $\set{P}^{(1)}, ..., \set{P}^{(N)}$ with each basis $\set{P}^{(n)} = \{\set{P}^{(n)}_0, \set{P}^{(n)}_1, ... \}$ being orthogonal w.r.t. $x_n$'s marginal PDF $dF_{n}$:

\begin{equation}
\int_{\Omega_n} \set{P}^{(n)}_i(x_n) \set{P}^{(n)}_j(x_n) \, dF_n(x_n) = 0 \,\, \forall n \mbox{ iff } i \ne j
\end{equation}

The PCE of bounded degree $D$ approximates a function $f: \Omega \subset \mathbb{R}^N \to \mathbb{R}$ as a truncated expansion in terms of these bases:

\begin{equation}
\label{eq:pce}
f(\vec{x}) \approx \sum_{\pmb{\alpha} = (0, ..., 0)}^{(D, ..., D)} \set{C}_{\pmb{\alpha}} \cdot \Psi_{\pmb{\alpha}}(\vec{x})
\end{equation}
with

\begin{equation}
\Psi_{\pmb{\alpha}}(\vec{x}) := \prod_{n=1}^N \set{P}^{(n)}_{\alpha_n}(x_n)
\end{equation}
being $\set{C}$ an $N$-dimensional tensor containing the expansion weights.

Sudret~\cite{Sudret:08} established a connection between the Sobol decomposition and the PCE that has gained significant popularity~\cite{GCAI:16, Bigoni:14, DKLM:14}. The author proposed the indices $SU_{\pmb{\alpha}}$, which approximate each Sobol coefficient $S_{\pmb{\alpha}}$ from a PCE surrogate of bounded degree. The idea behind $SU$, similar to~\cite{Sobol:69}, is based on the fact that the first function $\mathcal{P}^{(n)}_0$ from every PCE basis is a zero-degree polynomial and therefore a constant. From this and the basis' orthogonality it follows that the fraction of the model's response that is \emph{not} explained by a specific variable $n$ is exactly the one captured by the projection onto $\mathcal{P}^{(n)}_0$, while the remaining $\{\mathcal{P}^{(n)}_d\}_{d \ge 1}$ account for the interactions where the variable is present.

We can convert any such PCE representation into a TT surrogate as follows, at the expense of only a small discretization error (that can be easily adjusted, recall Sec.~\ref{sec:discretization}). Eq.~\ref{eq:pce} is interpretable as a \emph{continuous} Tucker decomposition, with the $f_{\pmb{\alpha}}$ acting as the elements of a core of size $(D+1)^N$. To obtain a standard Tucker format we just need to define its factor matrices. Each factor $\mat{U}^{(n)}$ has size $I_n \times (D+1)$ and is found by sampling the corresponding polynomial basis over the discretized variable range $x_n(1), ..., x_n(I_n)$:

\begin{equation}
\label{eq:tucker_discretization}
\mat{U}^{(n)}[i, j] := \set{P}^{(n)}_j(x_n(i))
\end{equation}
After this we can apply the conversion method detailed in Sec.~\ref{sec:from_tucker} to get the equivalent TT representation.

Alternatively, we may also convert a low-rank PCE expansion~\cite{KS:16} to TT by means of the CP conversion method above.
Such low-rank expansions define a \emph{continuous} CP as

\begin{equation}
f(\vec{x}) \approx \sum_{r=1}^R \lambda_r \cdot \Psi_r(\vec{x})
\end{equation}
with each rank-1 component being the outer product of functions that admit a low-degree polynomial expansion:

\begin{equation}
\Psi_r(\vec{x}) := \prod_{n=1}^N \left( \sum_{d=0}^D \alpha^{(n)}_{rd} \cdot \mathcal{P}^{(n)}_d(x_n) \right)
\end{equation}
and can be converted into standard PCE via discretization, analogously to Eq.~\ref{eq:tucker_discretization}.

\section{The Sobol Tensor Train} \label{sec:global}



We now introduce our proposed \emph{Sobol tensor train}, denoted as $\set{S}$, which has dimension $N$ and size 2 along each dimension for a total of $2^N$ elements. Such $2 \times ... \times 2$ tensors are not unusual, see for example the so-called \emph{quantized tensor train} (QTT) and the closely-related \emph{wavelet tensor train} (WTT)~\cite{OT:11}, as well as the recent multilinear regressors known as \emph{exponential machines}~\cite{NTO:16}. $\set{S}$ hence records the Sobol indices for all $n$-ary interactions:

\begin{equation}
S_{\pmb{\alpha}} \approx \set{S}_{\pmb{\alpha}} = \set{S}[j_1, \dots, j_N] = \set{S}^{(1)}[j_1] \cdot ... \cdot \set{S}^{(N)}[j_N]
\end{equation}
with $j_n = 1$ if $n \in \pmb{\alpha}$ and $0$ otherwise. 
To construct it we combine the definitions of Sobol decomposition (Sec.~\ref{sec:sobol_decomposition}) and Sobol indices (Sec.~\ref{sec:variance_and_sobol}) with the TT formulation as follows.

\begin{proposition} \label{prop:tt_sobol}
Let $\vec{x} = (x_1, ..., x_N)$ be independent with distributions $F_1, ..., F_N$, and let $\set{T} = [[\set{T}^{(1)}, ..., \set{T}^{(N)}]]$ be a TT surrogate $\tilde{f}(\vec{x}) \approx f(\vec{x})$. Then, each term $\tilde{f}_{\pmb{\alpha}}$ of the Sobol decomposition of $\tilde{f}$ is given by $\set{T}_{\pmb{\alpha}} = [[\set{T}^{(1)}_{\pmb{\alpha}}, ..., \set{T}^{(N)}_{\pmb{\alpha}}]]$ with cores defined slice-wise as

\begin{equation}
\label{eq:tt_sobol}
\set{T}^{(n)}_{\pmb{\alpha}}[j] :=
\begin{cases}
\operatorname{E}[\set{T}^{(n)}] \mbox{ if } n \notin \pmb{\alpha} \\
\set{T}^{(n)}[j] - \operatorname{E}[\set{T}^{(n)}] \mbox{ if } n \in \pmb{\alpha}
\end{cases}
\end{equation}
for all slices $j = 0,..., I_n-1$, where $\operatorname{E}[\set{T}^{(n)}] := \frac{1}{I_n} \cdot \sum_{i=0}^{I_n-1} F_n(x_n(i)) \, \set{T}^{(n)}[i]$ is the expectation along the $n$-th dimension, i.e. the average of the $n$-th core's slices, weighted by the $n$-th PDF term.
\end{proposition}

See the~\ref{app:proposition} for a proof. Eq.~\ref{eq:tt_sobol} can be intuitively interpreted as follows: Variables not in $\pmb{\alpha}$ must be integrated over their domain of existence, and $\tilde{f}_{\pmb{\alpha}}$ does not effectively depend on them. Their corresponding cores are accordingly constant. For variables in $\pmb{\alpha}$, on the other hand, we must keep the original function but subtract from it the lower-order expectations; these are all correctly accounted for thanks to multilinearity.

With Eq.~\ref{eq:tt_sobol} we can already compute any arbitrary Sobol index $\pmb{\alpha}$. However, we give now a more expeditious method that allows us to produce \emph{all} indices at once. First note that the following tensor $\set{T}_* = [[\set{T}^{(1)}_*, \dots, \set{T}^{(N)}_*]]$ compactly encodes all $2^N$ Sobol decomposition terms:

\begin{equation}
\label{eq:tt_sobol_all}
\begin{cases}
\set{T}^{(n)}_*[0] := \operatorname{E}[\set{T}^{(n)}] \\
\set{T}^{(n)}_*[j] := \set{T}^{(n)}[j-1] - \operatorname{E}[\set{T}^{(n)}] \mbox{ for } j = 1, ..., I_n+1
\end{cases}
\end{equation}

This tensor $\set{T}_*$ is simply a concatenation of the two types of slices of Eq.~\ref{eq:tt_sobol}, and so it approximates $f_{\pmb{\alpha}}(\vec{x})$ for any $\pmb{\alpha}$ and any $\vec{x} \in \Omega$.

We have now all necessary components to construct our Sobol tensor $\set{S}$: we need to compute $\operatorname{V}[\tilde{f}_{\pmb{\alpha}}]/D = \operatorname{E}[(\tilde{f}_{\pmb{\alpha}} - \operatorname{E}[\tilde{f}_{\pmb{\alpha}}])^2]/D$ in the compressed domain. The procedure, detailed in Alg.~\ref{alg:sobol}, also obtains the variance indices $D_{\pmb{\alpha}}$ as a necessary subproduct prior to normalization by the total variance $D$.

\begin{algorithm}[h]
\begin{algorithmic}[1]
\STATE Compute $\set{T}_*$ as in Eq.~\ref{eq:tt_sobol_all} \COMMENT{$\set{T}_*$ encodes $\tilde{f}_{\pmb{\alpha}} \, \forall \pmb{\alpha}$}
\STATE Compute $\set{T}_{**} := \set{T}_* \circ \set{T}_* = \set{T}_*^2$ \COMMENT{$\set{T}_{**}$ encodes $\tilde{f}^2_{\pmb{\alpha}} \, \forall \pmb{\alpha}$}
\FOR {$n = 1, \dots, N$}
	\STATE $\set{D}^{(n)}[0] := \set{T}_{**}^{(n)}[0]$
	\STATE $\set{D}^{(n)}[1] := \frac{1}{I_n} \cdot \sum_{i=0}^{I_n-1} F_n(x_n(i)) \, \set{T}_{**}^{(n)}[i+1]$
\ENDFOR
\STATE $\set{D} := [[\set{D}^{(1)}, ..., \set{D}^{(N)}]]$ \COMMENT{$\set{D}$ encodes $\operatorname{V}[f_{\pmb{\alpha}}] \, \forall \pmb{\alpha}$}
\STATE $D := \prod_{n=1}^N (\set{D}^{(n)}[0] + \set{D}^{(n)}[1])$ \COMMENT{Total variance $\operatorname{V}[\tilde{f}]$}
\STATE $\set{S} := \set{D} / D$
\RETURN $\set{S}$
\end{algorithmic}
\caption{Given a TT surrogate $\set{T} = [[\set{T}^{(1)}, ..., \set{T}^{(N)}]]$ of size $I_1 \times ... \times I_N$, extract the compressed Sobol tensor $\set{S}$ of size $2 \times ... \times 2$.}
\label{alg:sobol}
\end{algorithm}

If the input surrogate has TT-ranks $R_1, ..., R_{N-1}$, then $\set{S}$ may have at most ranks $R_1^2, ..., R_{N-1}^2$. The squaring (line 2 from Alg.~\ref{alg:sobol}) can be achieved either by ACA or by slice-wise Kronecker product if the rank is low enough (see~\ref{app:operations}). All other operations cannot increase any of the ranks. Last, note that the first corner coefficient in the tensor $\set{S}_{\emptyset} = \set{S}[0,\dots,0] = \tilde{f}_{\emptyset}/D$ is not a Sobol index; we set it to zero if needed with a simple rank-1 correction:

\begin{equation}
\set{S} \leftarrow \set{S} - \begin{pmatrix}\tilde{f}_{\emptyset}/D \\ 0\end{pmatrix} \otimes \overbrace{\begin{pmatrix}1 \\ 0\end{pmatrix} \otimes ... \otimes \begin{pmatrix}1 \\ 0\end{pmatrix}}^{N-1\mbox{ terms}}
\end{equation} 

%

\section{Computing Aggregated Sobol Indices} \label{sec:set_operations}

Aggregated indices require up to an exponential number of addends if computed naively. But thanks to multilinearity of the proposed tensor decomposition, we can obtain all such QoI at once and at very little cost.

\subsection{Superset Sobol Tensors} \label{sec:superset_sobol_tensor}

We recall now the notion of superset indices from Sec.~\ref{sec:related_indices}, which capture the aggregate dependence with respect to a group of indices; i.e. Sobol indices over $\pmb{\alpha}$ plus the indices of all variable tuples that are a superset of $\pmb{\alpha}$.
%
%
If $\set{S}$ is available, we can construct a superset Sobol tensor $\set{S}^S = [[\set{S}^{S (1)}, ..., \set{S}^{S (N)}]]$ that approximates any $S^S_{\pmb{\alpha}} \approx \set{S}^S_{\pmb{\alpha}}$. We construct its cores by slice-wise manipulation of the original cores:

\begin{equation}
\label{eq:s_to_ss}
\begin{cases}
\set{S}^{S (n)}[0] := \set{S}^{(n)}[0] + \set{S}^{(n)}[1] \\
\set{S}^{S (n)}[1] := \set{S}^{(n)}[1]
\end{cases}
\end{equation}

The rationale is that variables that are present in a tuple (encoded by the second slices, $j=1$) should stay present, while the rest (first slices, $j=0$) should be accounted for both when they are absent and when they are included. As an example, let us consider in 2D the second superset index $\set{S}^S_2 := \set{S}_{12} + \set{S}_2$. Eq.~\ref{eq:s_to_ss} yields

\begin{equation}
\begin{split}
\set{S}^S_{2} = \set{S}^{S (1)}[0] \cdot \set{S}^{S (2)}[1] = (\set{S}^{(1)}[0] + \set{S}^{(1)}[1]) \cdot \set{S}^{(2)}[1] \\ = \set{S}^{(1)}[0] \cdot \set{S}^{(2)}[1] + \set{S}^{(1)}[1] \cdot \set{S}^{(2)}[1] = \set{S}_1 + \set{S}_{12}
\end{split}
\end{equation}
as expected. Conversely, one may extract $\set{S}$ from $\set{S}^S$ by reverting the slice-wise transformations:

\begin{equation}
\label{eq:ss_to_s}
\begin{cases}
\set{S}^{(n)}[0] = \set{S}^{S (n)}[0] - \set{S}^{S (n)}[1] \\
\set{S}^{(n)}[1] = \set{S}^{S (n)}[1]
\end{cases}
\end{equation}

We wish to emphasize the compactness and convenience of the relations given by Eqs.~\ref{eq:s_to_ss} and~\ref{eq:ss_to_s}. A naive sum to obtain a superset index of order $K$ out of the standard indices $\set{S}$ would require $2^{N-K}$ addends. For example, for $N = 3$ and $\pmb{\alpha} = \{1\}$ we have $\set{S}^S_1 = \set{S}_1 + \set{S}_{12} + \set{S}_{13} + \set{S}_{123}$. Conversely, producing indices $\set{S}$ from $\set{S}^S$ needs $2^{N-K}$ mixed additions and subtractions as dictated by the inclusion-exclusion principle from combinatorics. For instance, $\set{S}_1 = \set{S}^S_1 - \set{S}^S_{12} - \set{S}^S_{13} + \set{S}^S_{123}$. On the other hand, Eqs.~\ref{eq:s_to_ss} and~\ref{eq:ss_to_s} need only $O(N R^2)$ operations in the TT format.

\subsection{Closed Sobol Tensors} \label{sec:inclusive_sobol_tensor}


Similarly to Eq.~\ref{eq:s_to_ss}, we derive the closed Sobol tensor $\set{S}^C$ from $\set{S}$ as follows:

\begin{equation}
\label{eq:s_to_sc}
\begin{cases}
\set{S}^{C (n)}[0] := \set{S}^{(n)}[0] \\
\set{S}^{C (n)}[1] := \set{S}^{(n)}[0] + \set{S}^{(n)}[1]
\end{cases}
\end{equation}

The logic here is that indices absent in a tuple should stay absent, while present indices should be accounted for also when they are missing (since we want to sum all subsets). The converse equation reads

\begin{equation}
\label{eq:sc_to_s}
\begin{cases}
\set{S}^{(n)}[0] = \set{S}^{C (n)}[0] \\
\set{S}^{(n)}[1] = \set{S}^{C (n)}[1] - \set{S}^{C (n)}[0]
\end{cases}
\end{equation}

\begin{figure}[t]\center
 \includegraphics[width=1\columnwidth]{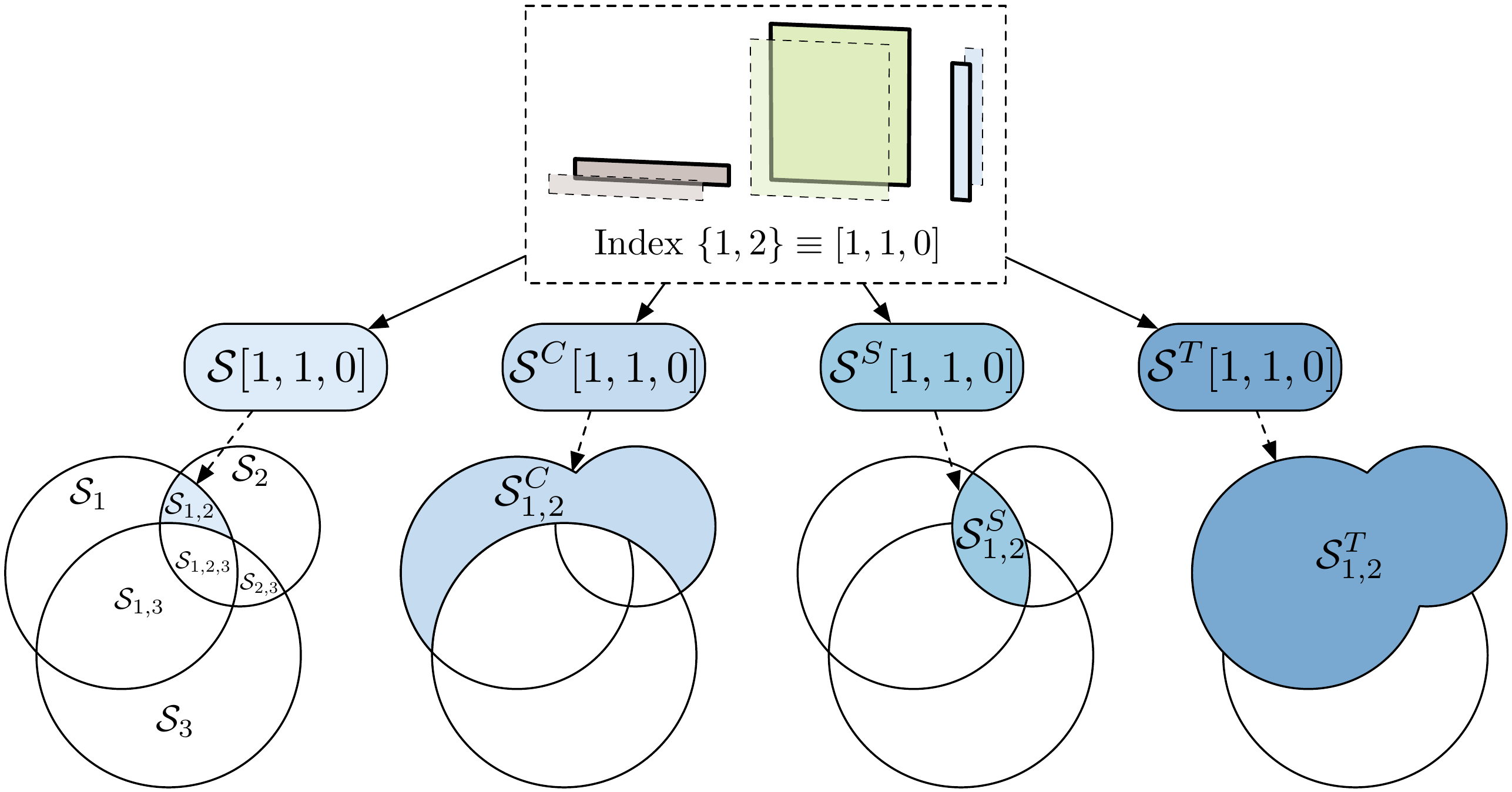}
 \caption{Examples of standard $\set{S}$, closed $\set{S}^C$, superset $\set{S}^S$, and total Sobol indices $\set{S}^T$ for a 3-variable model, interpreted as set cardinalities. Each colored region area is obtained from its corresponding tensor by multiplying together the indexed slices.}
  \label{fig:venn}
\end{figure}

\subsection{Total Sobol Tensors} \label{sec:total_sobol_tensor}

Our last aggregated tensor is the total $\set{S}^T$ and can be obtained via the complement operation as $\set{S}^T_{\pmb{\alpha}} = 1 - \set{S}^C_{-\pmb{\alpha}}$. Let us define a complement tensor $\bar{\set{S}}^C$, defined for each tuple as $\bar{\set{S}}^C_{\pmb{\alpha}} := \set{S}^C_{-\pmb{\alpha}}$. We extract this tensor from $\set{S}^C$ by simply swapping the two slices of each core:

\begin{equation}
\begin{cases}
\bar{\set{S}}^{C (n)}[0] := \set{S}^{C (n)}[1] \\
\bar{\set{S}}^{C (n)}[1] := \set{S}^{C (n)}[0]
\end{cases}
\end{equation}
and the final result $\set{S}^T = 1 - \bar{\set{S}}^C$ follows from a simple tensor-tensor subtraction. To retrieve $\set{S}^C$ back from $\set{S}^T$ it suffices to repeat the whole transformation.

\section{Global Sensitivity Metrics and Queries} \label{sec:other_indices}

\subsection{Relevant Subsets of Variables}

A typical and fundamental target in SA is to ``\emph{select the $k$ variables that account for the most variance}'', or alternatively ``\emph{select the smallest set variables that account for at least (say) 99\% variance}''. In order to tackle this we introduce the \emph{Hamming mask} of order $k$, that we define as

\begin{equation}
\set{M}^k_{\pmb{\alpha}} :=
\begin{cases}
1 & \mbox{ if } |\pmb{\alpha}| = k \\
0 & \mbox{ otherwise}
\end{cases}
\end{equation}



\begin{figure}[t]\center
	\includegraphics[width=1\columnwidth]{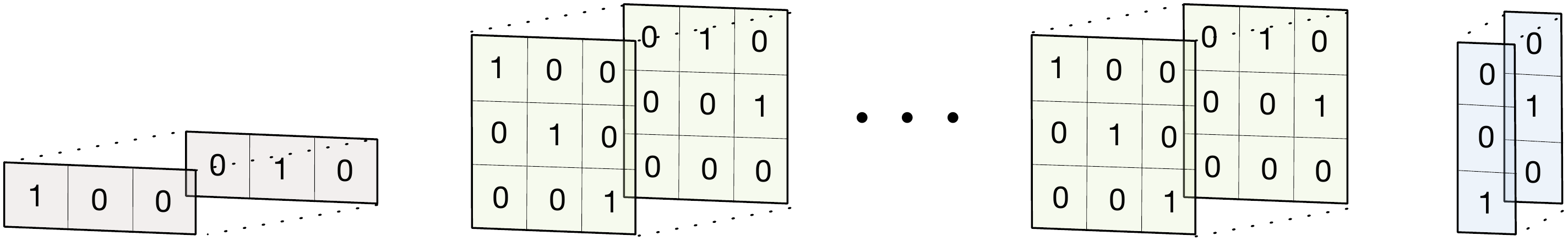}
	\caption{The Hamming mask tensor train $\set{M}^k$ for order $k=2$. At each position $\pmb{\alpha} \in \{0, 1\}^N$ it contains a `1' if and only if $|\pmb{\alpha}| = k$, and 0 otherwise. It is compressed with $N$ cores (rank $k+1$) using $2 (k+1)^2 (N-2) + 4(k+1)$ elements in total.}
	\label{fig:hamming_mask}
\end{figure}

We are able to build its compressed version using only $k+1$ ranks (Fig.~\ref{fig:hamming_mask}). It is best understood by following the vector-matrix sequence of products that takes place to reconstruct one element, left to right. The vector at the $n$-th step has size $k+1$ and encodes how many `1' bits have been encountered so far: a `1' at its first position means 0, a `0' followed by a `1' means 1, etc. The core slices transform this vector counter to account for the new bits as we traverse the binary sequence. The first slice of each core is the identity matrix, since it corresponds to a bit set to 0 (which does not have an effect). The second slice, however, must increment the counter, i.e. shift the `1' one position towards the right. It is therefore implemented as a shifted identity matrix. The last core simply checks if the total number of `1' found until the end matches $k$ or not.

The mask tensor $\set{M}_{\pmb{\alpha}}$ allows us to define restricted searches. For instance, the single most important Sobol index of order $k$ is

\begin{equation}
\argmax_{\pmb{\alpha}} \left\{ (\set{S} \circ \set{M}^k)_{\pmb{\alpha}} \right\}
\end{equation}
which we solve using a state-of-the-art global optimization algorithm in the TT format (Sec.~\ref{sec:global_optimization}). One may also use $\set{S}^T, \set{S}^C$ or $\set{S}^S$ instead of $\set{S}$ depending on the task at hand. For example, the $\set{S}_{\pmb{\alpha}}$ do play the dominant role in factor prioritization, but for factor fixing one is advised to seek a tuple with the smallest total index~\cite{SRACCGST:08}. 

We also use $\set{M}^k$ to compute the overall per-order contributions: the tensor dot product
\begin{equation} \label{eq:mask_dot}
<\set{S}, \set{M}^k>
\end{equation}
gives us the combined order $k$ indices $\sum_{|\pmb{\alpha}| = k} \set{S}_{\pmb{\alpha}}$.

\subsection{Other Constraints}

The analyst may seek a model simplification that satisfies additional constraints, e.g. that certain variables must, or must not, become frozen. Such conditions can be easily imposed by editing the mask tensor $\set{M}$. For instance, if a variable $1 \le n \le N$ should be fixed (i.e. simplified) it is sufficient to remove the second slice of the $n$-th mask core. This 
effectively restricts the search to $\{\pmb{\alpha} \given n \notin \pmb{\alpha}\}$ as desired. Conversely, if we wish to ensure that a variable is not fixed (i.e. remains active in the new simplified model), we just remove the first slice of the $n$-th core of $\set{M}$.

\section{Experimental Results} \label{sec:results}

Our experiments were conducted in Python. We exploit the ttpy toolbox\cite{ttpy}, a Python/Fortran library for TT manipulation that supports, among others, compression from full explicit tensors, slicing, decompression, truncation (rounding), and cross-approximation for any dimensionality.


\subsection{Sobol ``G'' Function}

\begin{figure}[h]\center
	\includegraphics[width=0.85\columnwidth]{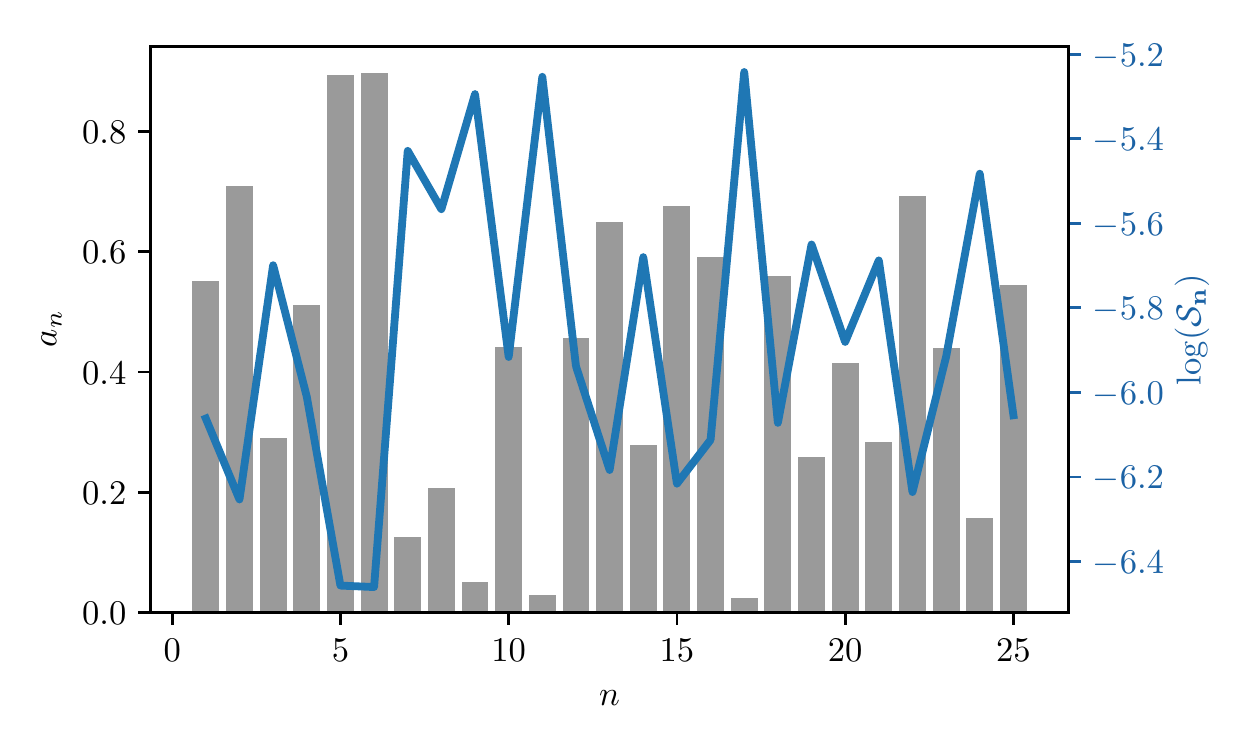}
	\caption{The 25 random values chosen for $a_n$ and their resulting Sobol values (depicted in logarithmic scale).}
	\label{fig:sobol_an}
\end{figure}

This function has been extensively used in the SA literature owing to its flexibility and relatively high-order interactions. It is defined as

\begin{equation}
\label{eq:sobol}
f(\vec{x}) := \prod_{n=1}^N \frac{|4 x_n - 2| + a_n}{1 + a_n}
\end{equation}
being $a_i$ random coefficients sampled from a uniform distribution $\mathcal{U}(0, 1)$ and $x_n \sim \mathcal{U}(0, 1)$ the  $n$ function parameters. Note that $f$ is non-differentiable at one point, namely $(0.5, ..., 0.5)$. We can expect to get an exact TT interpolator (up to machine precision) of $f$ as it is a product of univariate functions and therefore it has a rank 1. Our test example uses $N = 25$ dimensions and we discretize each variable into $I_1 = \dots = I_{25} = 64$ possible values. The ACA used 3200 evaluations of $f$ and was able to achieve a relative error of $\eps \approx 4.646 \cdot 10^{-15}$ over a test set of 4096 samples drawn at random using a Latin Hypercube sampling (LHS) scheme. Extracting the Sobol TT from the surrogate took $4.93$ seconds using $\eps = 10^{-6}$ as the ACA relative error for the squaring step in Alg.~\ref{alg:sobol}. Fig.~\ref{fig:sobol_an} shows every value of $a_n$ and its corresponding first-order Sobol index, computed using our method (Alg.~\ref{alg:sobol}).

Tab.~\ref{tab:Sobol1} shows the 5 highest Sobol indices of any order from our method, which in this case are only order-1 effects. We also show the indices as estimated directly from sampling the TT surrogate via the Sensitivity Analysis Library (SALib~\cite{HU:17}) in Python, using quasi-MC with varying number of sample points $P$. Finally and to complete the cross-check, we list the analytical Sobol values for comparison~\cite{Sobol:03}: $D_n = 1/(3(1+a_n)^2), D = \prod_n (D_n+1) - 1$, and $S_n = (\prod_n D_n)/D$. Tab.~\ref{tab:Sobol2} shows the highest aggregated indices (i.e. total, closed and superset) of order 1, 2, and 3 separately.

\begin{table*}[h]
	\centering
	\caption{Highest Sobol indices for the Sobol G function}

    \begin{tabular}{|c|c|cccc|}
        \hline
        \multicolumn{1}{|c|}{\multirow{3}{*}{\textbf{Index}}} & \multicolumn{5}{c|}{\textbf{Value}} \\ \cline{ 2-6 }
        \multicolumn{1}{|c|}{} & \multicolumn{1}{c|}{\multirow{2}{*}{Sobol TT}} & \multicolumn{3}{c|}{SALib} & \multicolumn{1}{c|}{{\multirow{ 2 }{*}{ Analytical }}} \\ \cline{ 3-5 }
        \multicolumn{1}{|c|}{} & \multicolumn{1}{c|}{} & \multicolumn{1}{c|}{($P= 520000$)} & \multicolumn{1}{c|}{($P= 5.2 \cdot 10^6$)} & \multicolumn{1}{c|}{($P= 5.2 \cdot 10^7$)} & \multicolumn{1}{c|}{} \\ \hline\hline
    
            $\mathcal{S}_{ 17 }$ & 0.0053 & 0.0147 & 0.0110 & 0.0073 & 0.0054 \\

            $\mathcal{S}_{ 11 }$ & 0.0052 & -0.0001 & 0.0053 & 0.0052 & 0.0054 \\

            $\mathcal{S}_{ 9 }$ & 0.0050 & 0.0237 & 0.0084 & 0.0053 & 0.0051 \\

            $\mathcal{S}_{ 7 }$ & 0.0044 & -0.0038 & 0.0035 & 0.0039 & 0.0045 \\

            $\mathcal{S}_{ 24 }$ & 0.0042 & 0.0259 & 0.0062 & 0.0042 & 0.0042 \\

        \hline
        \end{tabular}
	\label{tab:Sobol1}
\end{table*}

\begin{table*}[h]
	\centering
	\caption{Highest aggregated indices of order 1, 2, and 3 for the Sobol G function}

\begin{tabular}{|c|lll|}
   \hline
   \multicolumn{1}{|c|}{\multirow{2}{*}{\textbf{Order}}} & \multicolumn{3}{c|}{\textbf{Index / Variable(s)}} \\ \cline{2-4}
   \multicolumn{1}{|c|}{} & Total & Closed & Superset \\
   \hline \hline
    1 & $\mathcal{S}^T_{ 17 } = 0.2452$ & $\mathcal{S}^C_{ 17 }= 0.0053$ & $\mathcal{S}^S_{ 17 } = 0.2452$ \\ \hline

    2 & $\mathcal{S}^T_{ 11,17 } = 0.4296$ & $\mathcal{S}^C_{ 11,17 }= 0.0122$ & $\mathcal{S}^S_{ 11,17 } = 0.0586$ \\ \hline

    3 & $\mathcal{S}^T_{ 9,11,17 } = 0.5657$ & $\mathcal{S}^C_{ 9,11,17 }= 0.0209$ & $\mathcal{S}^S_{ 9,11,17 } = 0.0136$ \\  \hline
    \end{tabular}
	\label{tab:Sobol2}
\end{table*}

We observe that the TT indices are accurate to almost 4 decimal digits. The Sobol indices computed by SALib become closer the more samples are taken, further supporting the correctness of our method. Note that for this function SALib required a very large number of samples to obtain the Sobol indices with a similar level of precision, in contrast to our proposed method. We attribute this to the function's high dimensionality, which can be handled well by the TT.

\subsection{Piston Simulation}

This is a lower-dimensional but more complex model that measures the cycle time of a piston simulation~\cite{KZ:98}. The output is defined analytically on 7 variables as

\begin{equation}
f(\vec{x}) = 2 \pi \sqrt{\frac{M}{k + S^2 \frac{P_0 V_0}{T_0} \frac{T_a}{V^2}}}
\end{equation}
with

\begin{equation}
\begin{gathered}
V = \frac{S}{2k} \left(\sqrt{A^2 + 4k \frac{P_0 V_0}{T_0} T_a} - A\right) \\
A = P_0 S + 19.62 M - \frac{k V_0}{S}
\end{gathered}
\end{equation}

\begin{table*}[h]
	\centering
	\caption{Parameters of the Piston simulation}
	\begin{tabular}{c l c l}
		\hline
		\textbf{Variable} & \textbf{Description} & \textbf{Units} & \textbf{Distribution} \\
		\hline
		$M$ & Piston weight & ${kg}$ & $\mathcal{U}(30, 60)$ \\
		$S$ & Piston surface area & $m^2$ & $\mathcal{U}(0.005, 0.02)$ \\
		$V_0$ & Initial gas volume & $m^3$ & $\mathcal{U}(0.002, 0.01)$ \\
		$ k$ & Spring coefficient & $N/n$ & $\mathcal{U}(1000, 5000)$ \\
		$P_0$ & Atmospheric pressure & $N/m^2$ & $\mathcal{U}(90000, 110000)$ \\
		$T_a$ & Ambient temperature & $K$ & $\mathcal{U}(290, 296)$ \\
		$T_0$ & Filling gas temperature &$K$ & $\mathcal{U}(340, 360)$ \\
		\hline
	\end{tabular}
	\label{tab:Piston_params}
\end{table*}

The full list of parameters and their input ranges is detailed in Tab.~\ref{tab:Piston_params}. Our model was generated with ACA, stopped after 43904 function evaluations, again with $I = 64$ bins per dimension. It has 10496 non-zero elements and maximum rank $R = 7$, and it achieves $\eps \approx 0.077\%$ over an LHS-acquired test set. Note that the TT model is again built with fewer samples than those needed by SALib's MC algorithm, and that it is able to compute indices of arbitrary order \emph{a posteriori}. Extracting the Sobol TT took $5.70$ seconds in this case.

For further comparison we have also computed a PCE approximation of this function via 4096 training samples chosen similarly to the test set. To build the model we take the 4 first Legendre polynomials for each variable. We compress then the PCE-Tucker core into a TT model as detailed in Sec.~\ref{sec:from_low_rank} with a relative error of $0.5\%$, resulting in $R = 22$. The resulting TT-PCE model approximates the training set with a relative error $\eps \approx 0.38\%$, and achieves $\eps \approx 1.22\%$ on the test set.

As shown in Tab.~\ref{tab:Piston1}, our analysis reveals the fact that only the 4 first variables have a significant first-order effect. Their numerical values are consistent with the results reported e.g. in~\cite{ODC:14} (after normalization). Also, the most important tuple interactions arise from these very same variables. The triplet $\{S, V_0, k\}$ in particular has a closed index of about $95\%$ as reported in Tab.~\ref{tab:Piston2}. Overall, interactions of order 3 and above play a relatively small role. The proposed method is again able to compute all Sobol and aggregated indices in one go, and to do so with fewer evaluations than SALib.

\begin{table*}[h]
	\centering
	\caption{Highest Sobol indices for the piston function (interactions of order 3 and above are not supported by SALib)}

    \begin{tabular}{|c|c|cccc|}
        \hline
        \multicolumn{1}{|c|}{\multirow{3}{*}{\textbf{Index}}} & {\multirow{3}{*}{\textbf{Var(s)}}} & \multicolumn{4}{c|}{\textbf{Value}} \\ \cline{ 3-6 }
        \multicolumn{1}{|c|}{} &  & \multicolumn{1}{|c}{\multirow{2}{*}{Sobol TT}} & \multicolumn{1}{|c}{\multirow{2}{*}{Sobol TT-PCE}} & \multicolumn{1}{|c}{SALib on TT} & \multicolumn{1}{|c|}{{\multirow{ 1 }{*}{ SALib }}} \\ 
        \multicolumn{1}{|c|}{} &  & \multicolumn{1}{|c}{} & \multicolumn{1}{|c}{}  & \multicolumn{1}{|c}{ ($P$=160000) } & \multicolumn{1}{|c|}{ ($P$=160000) } \\ \hline\hline
   
            $\mathcal{S}_{ 2 }$ & $S$ & 0.5545 & 0.5585 & 0.5562 & 0.5563 \\

            $\mathcal{S}_{ 3 }$ & $V_0$ & 0.3207 & 0.3238 & 0.3215 & 0.3215 \\

            $\mathcal{S}_{ 1 }$ & $M$ & 0.0390 & 0.0396 & 0.0389 & 0.0391 \\

            $\mathcal{S}_{ 2,4 }$ & $S$,$k$ & 0.0242 & 0.0211 & 0.0252 & 0.0250 \\

            $\mathcal{S}_{ 4 }$ & $k$ & 0.0212 & 0.0200 & 0.0219 & 0.0221 \\

            $\mathcal{S}_{ 3,4 }$ & $V_0$,$k$ & 0.0129 & 0.0117 & 0.0121 & 0.0118 \\

            $\mathcal{S}_{ 2,3,4 }$ & $S$,$V_0$,$k$ & 0.0094 & 0.0066 & - & - \\

            $\mathcal{S}_{ 1,3 }$ & $M$,$V_0$ & 0.0050 & 0.0046 & 0.0053 & 0.0053 \\

            $\mathcal{S}_{ 2,3 }$ & $S$,$V_0$ & 0.0046 & 0.0043 & 0.0045 & 0.0044 \\

            $\mathcal{S}_{ 1,2 }$ & $M$,$S$ & 0.0046 & 0.0048 & 0.0036 & 0.0035 \\

        \hline
        \end{tabular}
	\label{tab:Piston1}
\end{table*}

\begin{table*}[h]
	\centering
	\caption{Highest aggregated indices of order 1, 2, and 3 for the piston function}

\begin{tabular}{|c|lll|}
   \hline
   \multicolumn{1}{|c|}{\multirow{2}{*}{\textbf{Order}}} & \multicolumn{3}{c|}{\textbf{Index / Variable(s)}} \\ \cline{2-4}
   \multicolumn{1}{|c|}{} & Total & Closed & Superset \\
   \hline \hline
    1 & $\mathcal{S}^T_{ 2 } = 0.5987$ & $\mathcal{S}^C_{ 2 }= 0.5545$ & $\mathcal{S}^S_{ 2 } = 0.5987$ \\ 
            & $\{$$S$$\}$ & $\{$$S$$\}$ & $\{$$S$$\}$ \\ \hline

    2 & $\mathcal{S}^T_{ 2,3 } = 0.9374$ & $\mathcal{S}^C_{ 2,3 }= 0.8799$ & $\mathcal{S}^S_{ 2,4 } = 0.0343$ \\ 
            & $\{$$S$,$V_0$$\}$ & $\{$$S$,$V_0$$\}$ & $\{$$S$,$k$$\}$ \\ \hline

    3 & $\mathcal{S}^T_{ 1,2,3 } = 0.9776$ & $\mathcal{S}^C_{ 2,3,4 }= 0.9475$ & $\mathcal{S}^S_{ 2,3,4 } = 0.0098$ \\ 
            & $\{$$M$,$S$,$V_0$$\}$ & $\{$$S$,$V_0$,$k$$\}$ & $\{$$S$,$V_0$,$k$$\}$ \\ \hline
    \end{tabular}
	\label{tab:Piston2}
\end{table*}

\subsection{GEMM Matrix Product in the GPU} \label{sec:gemm}

Our last experiment is a parallel computing example: we measured the computation time of 32-bit floating point matrix-matrix products in a graphics processing unit (GPU) according to 14 parameters and optimization techniques (loop unrolling, thread block-size, vector data types, etc.). The input variables are essentially discrete, since they are highly non-linear~\cite{NC:15} and can only take a handful of different values at most (usually a few powers of 2). We have chosen to build our TT surrogate using ALS tensor completion as we described in Sec.~\ref{sec:categorical_data}. The product analyzed is $\mat{A} \cdot \mat{B} = \mat{C}$ with all three matrices having size $2048 \times 2048$. We use the highly-tuneable GEMM kernel provided in the package CLTune~\cite{NC:15}, a generic auto-tuner for OpenCL kernels written in C++. Tab.~\ref{tab:CLTune_params} summarizes the 14 parameters and their input ranges (see the CLTune paper for further details). 

\begin{table*}[h]
	\centering
	\caption{The 14 parameters of the GEMM OpenCL Kernel}
	\begin{tabular}{c l c}
		\textbf{Variable(s)} & \textbf{Description} & \textbf{Domain} \\
		\hline
		$M_{wg}, N_{wg}$ & \specialcell{Per-matrix 2D tiling at \\workgroup level} & $\{16, 32, 64, 128\}$ \\ \hline
		$K_{wg}$ & \specialcell{Inner dimension of 2D tiling \\ at workgroup level} & $\{16, 32\}$ \\ \hline
		$M_{dimC}, N_{dimC}$ & Local workgroup size & $\{8, 16, 32\}$ \\ \hline
		$M_{dimA}, N_{dimB}$ & \specialcell{Local memory shape (when enable)} & $\{8, 16, 32\}$ \\ \hline
		$K_{wi}$ & Kernel loop unrolling factor & $\{2, 8\}$ \\	\hline
		$M_{vec}, N_{vec}$ & \specialcell{Per-matrix vector widths \\for loading and storing} & $\{1, 2, 4, 8\}$ \\ \hline
		$M_{stride}, N_{stride}$ & \specialcell{Enable stride for accessing off-chip \\memory within a single thread} & $\{yes, no\}$ \\ \hline
		$L\$_A, L\$_B$ & \specialcell{Per-matrix manual caching \\of the 2D workgroup tile} & $\{yes, no\}$ \\ \hline
	\end{tabular}
	\label{tab:CLTune_params}
\end{table*}

\begin{figure}[h]\center
	\includegraphics[width=0.9\columnwidth]{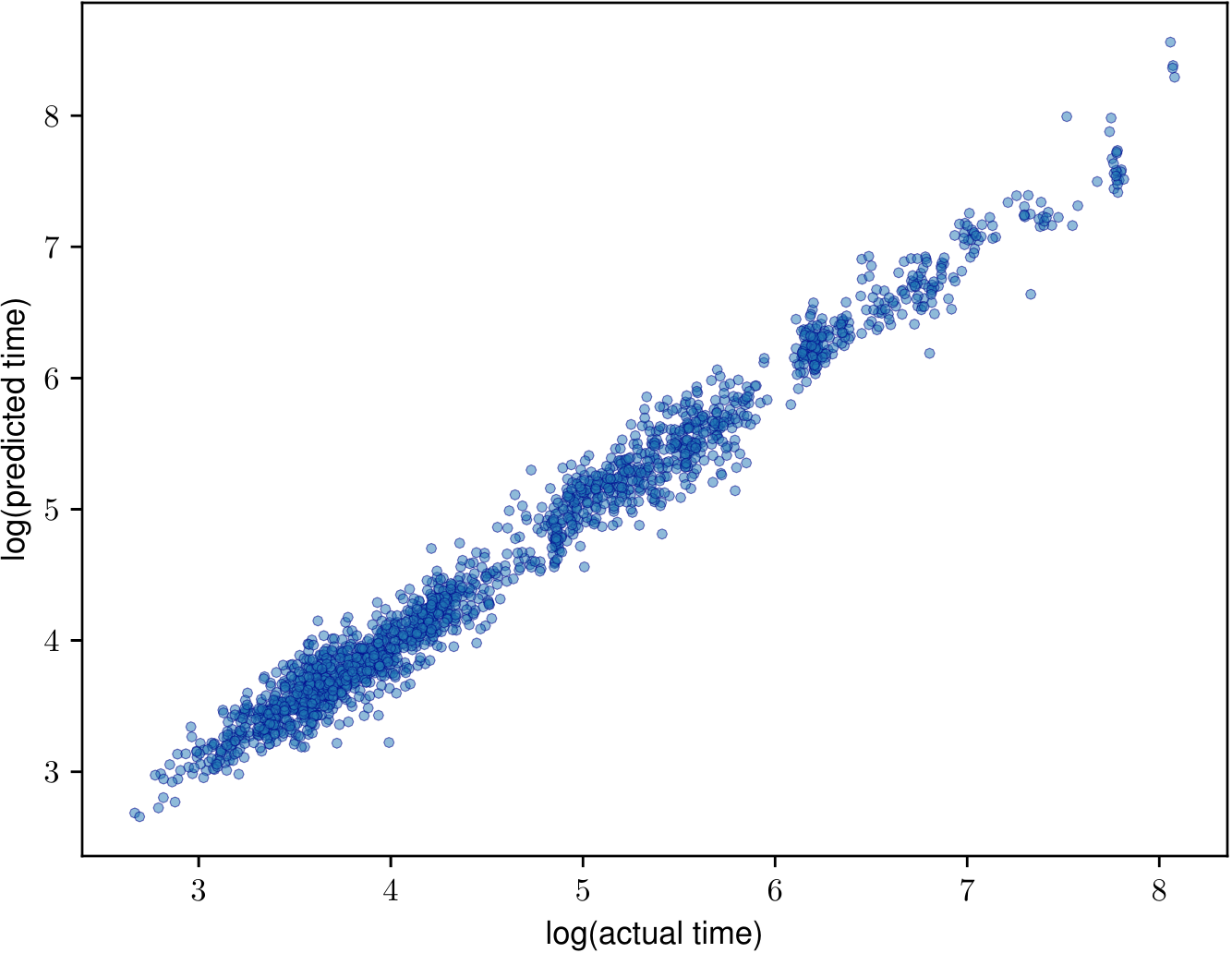}
	\caption{Surrogate obtained via TT completion for our GEMM experiment: groundtruth vs. prediction over the test set (1812 points), with relative error $\eps \approx 3.4\%$}
	\label{fig:gemm_plot_log}
\end{figure}

We generated this data set with a workstation running Ubuntu Linux 16.04, equipped with an Intel Core i5-4690 3.5GHz processor and a GeForce GTX680 GPU with 4GB of memory. The data set consists of 12080 samples taken uniformly at random (without repetition) among the 1327104 total possible variable combinations. Each sample was measured 25 times and averaged in order to reduce noise effects. All GEMM running times are considered in logarithmic scale both for training and analysis as advised in~\cite{FE:17}. We split the data as $70\%$, $15\%$, and $15\%$ for training, validation and test, respectively. The best TT surrogate was obtained after 25 ALS iterations, which took $28.7$ seconds. It has ranks $R_1 = \dots = R_{13} = 8$ for a total of 2224 non-zero elements. It achieved a relative error of $\epsilon \approx 3.4\%$ on the test set (see Fig.~\ref{fig:gemm_plot_log}). The final Sobol indices (Tables~\ref{tab:CLTune1} and \ref{tab:CLTune2}) were computed after re-fitting this best model to the full data set. The Sobol tensor took $12.32$ seconds to build.

\begin{table*}[h]
	\centering
	\caption{Highest Sobol indices for the GEMM matrix product function}

\begin{tabular}{|c|c|cc|}
\hline
\multicolumn{1}{|c|}{\multirow{3}{*}{\textbf{Index}}} & {\multirow{3}{*}{\textbf{Var(s)}}} & \multicolumn{2}{c|}{\textbf{Value}} \\ \cline{ 3-4 }
\multicolumn{1}{|c|}{} &  & \multicolumn{1}{c}{\multirow{2}{*}{Sobol TT}} & \multicolumn{1}{|c|}{SALib on TT} \\
\multicolumn{1}{|c|}{} &  & \multicolumn{1}{c}{} & \multicolumn{1}{|c|}{ ($P$=300000) } \\ \hline\hline

    $\mathcal{S}_{ 1,2,4 }$ & $M_{wg}$,$N_{wg}$,$M_{dimC}$ & 0.0842 & - \\

    $\mathcal{S}_{ 1,4 }$ & $M_{wg}$,$M_{dimC}$ & 0.0790 & 0.0799 \\

    $\mathcal{S}_{ 2,4 }$ & $N_{wg}$,$M_{dimC}$ & 0.0675 & 0.0754 \\

    $\mathcal{S}_{ 1 }$ & $M_{wg}$ & 0.0643 & 0.0539 \\

    $\mathcal{S}_{ 1,2 }$ & $M_{wg}$,$N_{wg}$ & 0.0628 & 0.0686 \\

    $\mathcal{S}_{ 1,4,5 }$ & $M_{wg}$,$M_{dimC}$,$N_{dimC}$ & 0.0330 & - \\

    $\mathcal{S}_{ 1,2,4,5 }$ & $M_{wg}$,$N_{wg}$,$M_{dimC}$,$N_{dimC}$ & 0.0319 & - \\

    $\mathcal{S}_{ 4 }$ & $M_{dimC}$ & 0.0286 & 0.0257 \\

    $\mathcal{S}_{ 4,5 }$ & $M_{dimC}$,$N_{dimC}$ & 0.0225 & 0.0165 \\

    $\mathcal{S}_{ 2 }$ & $N_{wg}$ & 0.0198 & 0.0051 \\
\hline
\end{tabular}
	\label{tab:CLTune1}
\end{table*}

\begin{table*}[h]
	\centering
	\caption{Highest aggregated indices of order 1, 2, and 3 for the GEMM matrix product}

\begin{tabular}{|c|lll|}
   \hline
   \multicolumn{1}{|c|}{\multirow{2}{*}{\textbf{Order}}} & \multicolumn{3}{c|}{\textbf{Index / Variable(s)}} \\ \cline{2-4}
   \multicolumn{1}{|c|}{} & Total & Closed & Superset \\
   \hline \hline
    1 & $\mathcal{S}^T_{ 1 } = 0.6979$ & $\mathcal{S}^C_{ 1 }= 0.0643$ & $\mathcal{S}^S_{ 1 } = 0.6979$ \\ 
            & $\{$$M_{wg}$$\}$ & $\{$$M_{wg}$$\}$ & $\{$$M_{wg}$$\}$ \\ \hline

    2 & $\mathcal{S}^T_{ 1,4 } = 0.8960$ & $\mathcal{S}^C_{ 1,4 }= 0.1718$ & $\mathcal{S}^S_{ 1,4 } = 0.4380$ \\ 
            & $\{$$M_{wg}$,$M_{dimC}$$\}$ & $\{$$M_{wg}$,$M_{dimC}$$\}$ & $\{$$M_{wg}$,$M_{dimC}$$\}$ \\ \hline

    3 & $\mathcal{S}^T_{ 1,2,4 } = 0.9540$ & $\mathcal{S}^C_{ 1,2,4 }= 0.4060$ & $\mathcal{S}^S_{ 1,2,4 } = 0.2301$ \\ 
            & $\{$$M_{wg}$,$N_{wg}$,$M_{dimC}$$\}$ & $\{$$M_{wg}$,$N_{wg}$,$M_{dimC}$$\}$ & $\{$$M_{wg}$,$N_{wg}$,$M_{dimC}$$\}$ \\ \hline
    \end{tabular}
	\label{tab:CLTune2}
\end{table*}

Our results indicate a relatively large presence of high-order interactions; this matches the prior knowledge that GPU kernel optimization is a challenging high-dimensional parameter space, and that the parameters' influences tend to be highly inter-dependent~\cite{NC:15}. In particular the most important order-1 Sobol index (from $M_{wg}$) is only about $6\%$, and all order-1 indices combined explain only less than one fourth of the total model variability. We also use this real-world data set to test our querying routines; we report some sample results in Tab.~\ref{tab:queries} involving various aggregated indices.

\begin{table*}[h]
	\centering
	\caption{Once the Sobol TT is available, we can satisfy efficiently various types of queries as detailed in Sec.~\ref{sec:other_indices} using mask tensors, constrained search and TT global optimization}
	\begin{tabular}{l c c c}
		\textbf{Query} & \textbf{Result} & \textbf{Value} & \textbf{Computing time (s)}\\
		\hline
		\specialcell{Variable that interacts \\ the most with $\{L\$_A, L\$_B\}$} & $M_{wg}$ & \specialcell{$\mathcal{S}^S_{1,13,14} = 0.0263$ \\ (as high as possible)} & $0.3598$ \\ \hline
		\specialcell{Variable that interacts \\ the least with $M_{wg}$} & $N_{stride}$ & \specialcell{$\mathcal{S}^S_{1,12} = 0.0069$ \\ (as low as possible)} & $0.4237$ \\ \hline
		\specialcell{Highest closed 3-tuple \\ that avoids $M_{wg}$} & $N_{wg}$, $M_{dimC}$, $N_{dimC}$ & \specialcell{$\mathcal{S}^C_{2,4,5} = 0.1828$ \\ (as high as possible)} & $0.2786$
		\\ \hline
		\specialcell{Highest closed 3-tuple \\ that includes $M_{wg}$} & $M_{wg}$, $N_{wg}$, $M_{dimC}$ & \specialcell{$\mathcal{S}^C_{1,2,4} = 0.5940$ \\ (as high as possible)} & $0.3686$
		\\ \hline
		\specialcell{6 variables that \\ can be frozen with \\ the least impact} & \specialcell{$K_{wg}, K_{wi}, M_{vec},$ \\ $N_{vec}, M_{stride}, N_{stride}$} & \specialcell{$\mathcal{S}^T_{3,8,9,10,11,12} = 0.2365$ \\ (as low as possible)} & $0.7307$ \\ \hline
	\end{tabular}
	\label{tab:queries}
\end{table*}

\begin{figure}[h!]\center
	\subfigure[Sobol G function]{\includegraphics[width=0.49\columnwidth]{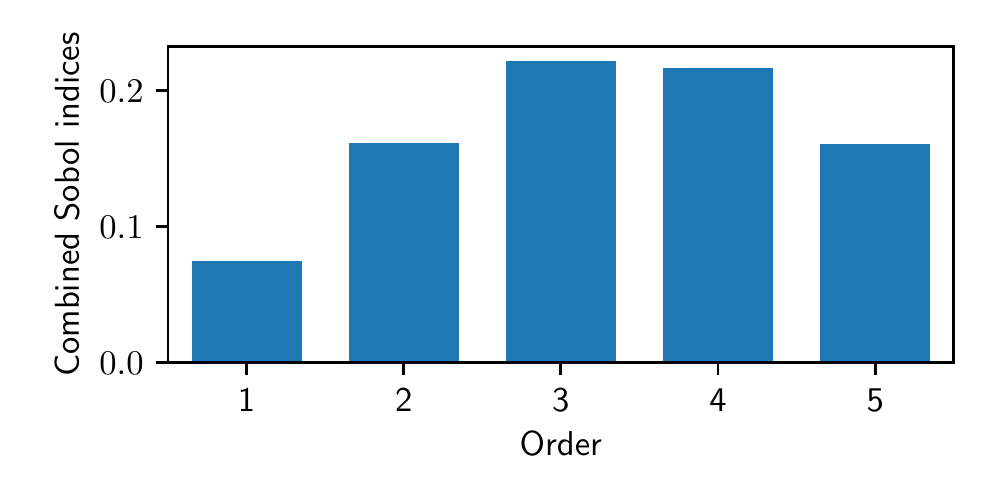}}
	\subfigure[Piston]{\includegraphics[width=0.49\columnwidth]{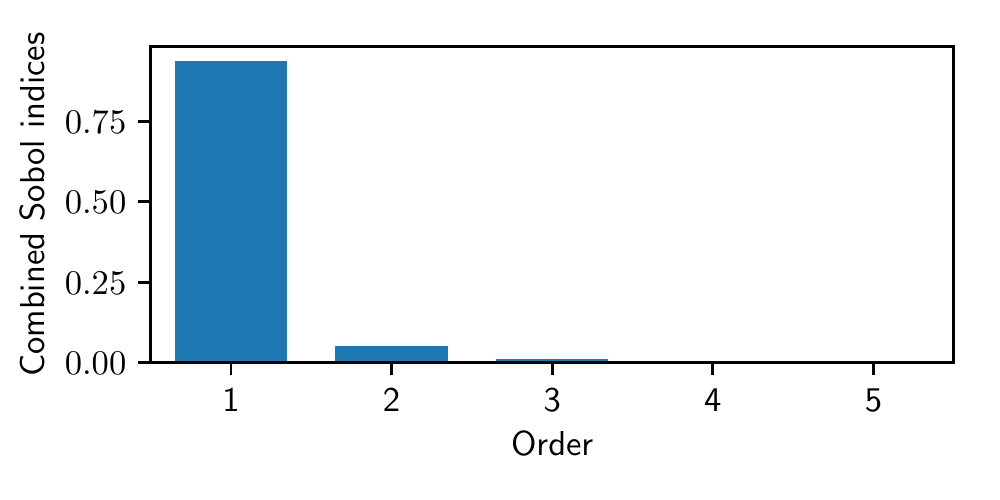}}
	\subfigure[GEMM]{\includegraphics[width=0.49\columnwidth]{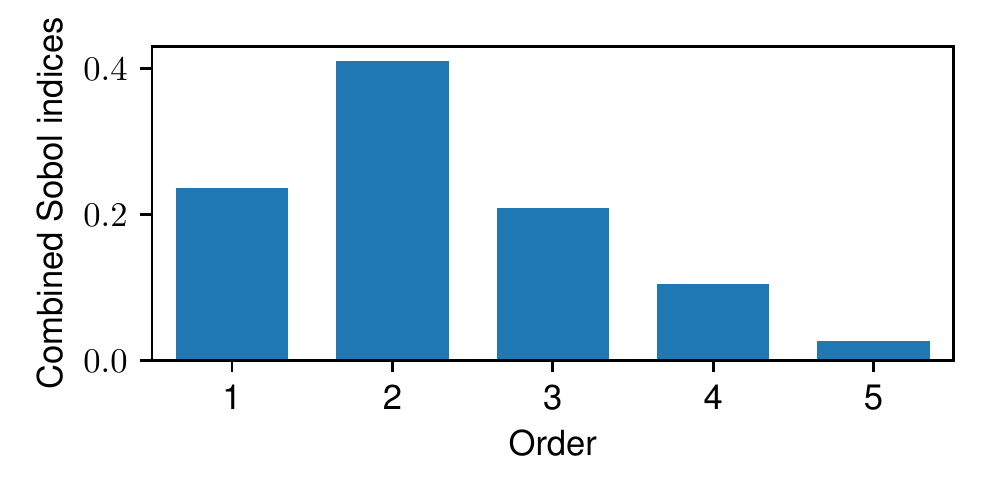}}
	\caption{Combined contributions for orders 1 to 5 for all three models, efficiently computed using Eq.~\ref{eq:mask_dot}}
	\label{fig:order_histograms}
\end{figure}

To conclude this section we show in Fig.~\ref{fig:order_histograms} one bar chart per data set, containing the overall relative variance broken down by interaction order.

\section{Conclusions}

We have introduced a compact data structure that gathers all Sobol indices from any TT-based surrogate model, and have given algorithms to extract various aggregated indices from it. The proposed aggregation algorithms capitalize on the format's multilinearity and have very little overhead cost. We combine these ideas with mask tensors, which allow us to define restricted queries and thus aid in model reduction/interpretation tasks. We believe the tensor train has a great potential as a canonical format for approximation of multiparametric systems, and the proposed methods for sensitivity analysis can be understood in the context of this trend. The presented framework is flexible in a variety of settings, and supports arbitrary orders of significant variable interactions in higher-dimensional models.

\subsection{Future Work}

Several possible extensions remain to be explored. For example, higher-order Sobol indices~\cite{ODC:14, GCAI:16} are useful for analysis of extreme values and risk minimization, and could be in principle ported to the TT format. Also, we wish to investigate more deeply TT surrogates for multi-valued models. Rather than training a separate model per individual output, we would like to work on a single tensor with an extra dimension to index the outputs. We believe that one may then run a joint TT analysis on the Sobol indices for all outputs at once and thus aid in model interpretability.


\appendix
\renewcommand{\thesection}{Appendix \Alph{section}}

\section{- Operations in the TT Format} \label{app:operations}

Multiplication/division of a TT tensor by a scalar $\alpha$ is achieved by simply multiplying/dividing one of its cores (say, the first) by $\alpha$. Tensor-tensor addition is written as $(\set{T}_1 + \set{T}_2)[\vec{x}] := \set{T}_1[\vec{x}] + \set{T}_2[\vec{x}]$ and has the following cores:

\begin{equation*}
\begin{cases}
\mleft(
\begin{array}{c|c}
  \set{T}_1^{(1)} [x_1] & \set{T}_2^{(1)} [x_1] \\
\end{array}
\mright) & \mbox{(first core)} \vspace{2mm} \\
\mleft(
\begin{array}{c|c}
  \set{T}_1^{(n)} [x_n] & 0 \\
  \hline
  0 & \set{T}_2^{(n)} [x_n]
\end{array}
\mright) & (1 < n < N) \vspace{2mm} \\
\mleft(
\begin{array}{c}
  \set{T}_1^{(N)} [x_N] \\
  \hline
  \set{T}_2^{(N)} [x_N]
\end{array}
\mright) & \mbox{(last core)} \vspace{2mm} \\
\end{cases}
\end{equation*}

The element-wise (or Hadamard) product $(\set{T}_1 \circ \set{T}_2)[\vec{x}] := \set{T}_1[\vec{x}] \cdot \set{T}_2[\vec{x}]$ arises from a slice-wise Kronecker product:

\begin{equation}
(\set{T}_1^{(1)}[x_1] \otimes \set{T}_2^{(1)}[x_1]) \cdot ... \cdot (\set{T}_1^{(N)}[x_N] \otimes \set{T}_2^{(N)}[x_N])
\end{equation}

\section{- Proof of Proposition~\ref{prop:tt_sobol}} \label{app:proposition}

Consider a tuple $\pmb{\alpha}$ and an arbitrary sampling point $\vec{x} = (x_1(i_1), \dots, x_N(i_N))$. We have defined our TT approximation as

\begin{equation} \label{eq:proof1}
\set{T}_{\pmb{\alpha}}(\vec{x}) = \prod_{n=1}^N \set{T}^{(n)}_{\pmb{\alpha}}[i_n]
\end{equation}
with

\begin{equation} \label{eq:proof2}
\set{T}^{(n)}_{\pmb{\alpha}}[i_n] :=
\begin{cases}
\operatorname{E}[\set{T}^{(n)}] \mbox{ if } n \notin \pmb{\alpha} \\
\set{T}^{(n)}[i_n] - \operatorname{E}[\set{T}^{(n)}] \mbox{ if } n \in \pmb{\alpha}
\end{cases}
\end{equation}

Expanding the $\pmb{\alpha}$ subtractions from Eq.~\ref{eq:proof2} we get a sequence of $2^{|\pmb{\alpha}|}$ additions and subtractions:

\begin{equation} \label{eq:proof3}
\sum_{\pmb{\beta} \subseteq \pmb{\alpha}} (-1)^{|\pmb{\alpha}| - |\pmb{\beta}|} \prod_{n=1}^N \widehat{\set{T}}^{(n)}_{\pmb{\beta}}[i_n]
\end{equation}
with

\begin{equation}
\widehat{\set{T}}_{\pmb{\beta}}^{(n)}[i_n] :=
\begin{cases}
\operatorname{E}[\set{T}^{(n)}] \mbox{ if } n \notin \pmb{\beta} \\
\set{T}^{(n)}[i_n] \mbox{ if } n \in \pmb{\beta}
\end{cases}
\end{equation}

Recall that $\set{T}^{(n)}$ encodes the model $\tilde{f}$'s response along the $n$-th axis, while $\operatorname{E}[\set{T}^{(n)}]$ represents its integration along that axis. Therefore Eq.~\ref{eq:proof3} becomes

\begin{equation}
\begin{split}
\sum_{\pmb{\beta} \subseteq \pmb{\alpha}} (-1)^{|\pmb{\alpha}| - |\pmb{\beta}|} \int_{\Omega - \pmb{\beta}} \tilde{f}(\vec{x}) dF_{-\pmb{\alpha}}(\vec{x}_{-\pmb{\alpha}}) \\
= \int_{\Omega-{\pmb{\alpha}}} \tilde{f}(\vec{x}) \, dF_{-\pmb{\alpha}}(\vec{x}_{-\pmb{\alpha}}) - \sum_{\pmb{\beta} | \pmb{\beta} \subset \pmb{\alpha}} \tilde{f}_{\pmb{\beta}}(\vec{x}_{\pmb{\beta}}) \\
= \tilde{f}_{\pmb{\alpha}}(\vec{x}_{\pmb{\alpha}}) \, \square
\end{split}
\end{equation}

\bibliographystyle{spphys}       

\bibliography{references}

\end{document}